\documentstyle[12pt, twoside]{article}
\input amssym.def
\input amssym.tex

\newenvironment{demo}[1]%
{\vskip-\lastskip\medskip
  \noindent
  {\em #1.}\enspace
  }%
{\qed\par\medskip
  }

\newcommand{\qed}{
  \strut\hfill
  \mbox{$\Box$}
  }

\newcommand{\ale}{ \widetilde{ \C^2 / \G} }
\newcommand{\alen}{ \ale^{[n]} }
\newcommand{\C}{ {\Bbb C} }
\newcommand{\End}{ \mbox{End} }
\newcommand{\Z}{ {\Bbb Z} }

\newcommand{\g}{\gamma}
\newcommand{\G}{\Gamma}
\newcommand{\Gn}{\G_n}

\newcommand{\hilq}{X // \G}
\newcommand{\hilqg}{ \C^{2} // \G }
\newcommand{\hilqgn}{ \C^{2n} // \Gn }
\newcommand{\hilqn}{ \C^{2n} // S_n }
\newcommand{\hnew}{X_{\G,n}}
\newcommand{\Hom}{ \mbox{Hom} }

\newcommand{\Ind}{\mbox{Ind} }
\newcommand{\Kgx}{ {K}^{\topo}_{\G}(X) }
\newcommand{\Kgxn}{ K^{\topo}_{\Gn}(X^n) }
\newcommand{\Ky}{ K(Y) }
\newcommand{\Kyn}{ K_{S_n} (Y^n) }
\newcommand{\la}{\lambda}
\newcommand{\affineg}{\widehat{\frak g}}
\newcommand{\quotuniv}{{\cal U}_{\G, n} }
\newcommand{\RG}{R_{\G} }
\newcommand{\Rz}{R_{\Bbb Z} (\G)}

\newcommand{\topo}{ {\footnotesize {top}} }
\newcommand{\univ}{ U_{\G,n} }
\newcommand{\x}{ {\bf x} }
\newcommand{\y}{ {\bf y} }
\newcommand{\wt}{ \xi }

\newtheorem{theorem}{Theorem}
\newtheorem{lemma}{Lemma}
\newtheorem{remark}{Remark}
\newtheorem{conjecture}{Conjecture}

\newtheorem{proposition}{Proposition}
\newtheorem{corollary}{Corollary}

\begin{document}
\title{Hilbert schemes, wreath products, and the McKay correspondence}
\author{
  Weiqiang Wang
\thanks{1991 {\em Mathematics Subject Classification}.
Primary 19J, 18F.
        }
}
\date{}
\maketitle

\begin{abstract}{
Various algebraic structures have recently appeared in a parallel
way in the framework of Hilbert schemes of points on a surface and
respectively in the framework of equivariant K-theory \cite{N1,
Gr, S2, W}, but direct connections are yet to be clarified to
explain such a coincidence. We provide several non-trivial steps
toward establishing our main conjecture on the isomorphism between
the Hilbert quotient of the affine space $\C^{2n}$ by the wreath
product $\Gn = \G \sim S_n$ and Hilbert schemes of points on the
minimal resolution of a simple singularity $\C^2 /\G$. We discuss
further various implications of our main conjecture. We obtain a
key ingredient toward a direct isomorphism between two forms of
McKay correspondence in terms of Hilbert schemes \cite{N1, Gr, N2}
and respectively of wreath products \cite{FJW}. We in addition
establish a direct identification of various algebraic structures
appearing in two different setups of equivariant K-theory
\cite{S2, W}. }
\end{abstract}
%\tableofcontents
%\setcounter{section}{-1}
%%
%%
%%
%%
%%
%%
%%
\section*{Introduction}
Nakajima \cite{N1} constructed a Heisenberg algebra using
correspondence varieties which acts on the direct sum over all $n$
of homology groups $H(X^{[n]})$ with complex coefficient of
Hilbert schemes $X^{[n]}$ of $n$ points on a quasi-projective
surface $X$. This representation is irreducible thanks to
G\"ottsche's earlier work \cite{G}. Similar results have been
independently obtained by Grojnowski \cite{Gr}. We refer to
\cite{N2} for an excellent account of Hilbert schemes and related
works. In the case when $X$ is the minimal resolution $\ale$ of
the simple singularity $\Bbb C^2 / \G$ associated to a finite
subgroup $\G$ of $SL_2 ( \C)$, this together with some additional
simple data provides a geometric realization of the
Frenkel-Kac-Segal vertex construction of the basic representation
of an affine Lie algebra \cite{FK, S1}. We may view this as a
geometric McKay correspondence \cite{Mc} which provides a
bijection between finite subgroups of $SL_2 (\C)$ and affine Lie
algebras of ADE types.

In \cite{W} we realized the important role of wreath products $\Gn
= \G \sim S_n $ associated to a finite group $\G$ in equivariant
K-theory. Various algebraic structures were constructed on the
direct sum over all $n$ of the topological $\Gn$-equivariant
K-theory $K^{\topo}_{\Gn} (X^n) \bigotimes \C$ of $X^n$ for a
$\G$-space $X$. The results of \cite{W} generalized the work of
Segal \cite{S2} (also see \cite{VW, Gr}) which corresponds to our
special case when $\G$ is trivial (i.e. $\G$ is the one-element
group) and the wreath product $\Gn$ reduces to the symmetric group
$S_n$.

The wreath product approach obtains further significance in light
of the conjectural equivalence of various algebraic structures in
the following three spaces:
\begin{eqnarray}   \label{eq_master}
 \begin{array}{ccc}
   \mbox{I} &  & \mbox{II}    \\
            \bigoplus_{n \geq 0} H( Y^{[n]})
          & \leftarrow  -\rightarrow
          &  \bigoplus_{n \geq 0} K^{\topo}_{S_n } (Y^n) \bigotimes \C \\
           & \;\, \nwarrow \quad \quad &\uparrow  \\
           & \quad \searrow & \downarrow  \\
          & &  \mbox{III}  \\
          & & \bigoplus_{n \geq 0} K^{\topo}_{\Gn } (X^n) \bigotimes \C
 \end{array}
\end{eqnarray}
Here one assumes that $X$ is a quasi-projective surface acted upon
by a finite group $\G$ and $Y$ is a suitable resolution of
singularities of $X/ \G$ such that there exists a canonical
isomorphism between $K_{\G}(X)$ and $K( Y)$. For $\G$ trivial
III reduces to II.
The graded dimensions of the three spaces have been shown to
coincide \cite{W}. The complexity of the geometry involved
decreases significantly from I to II, and then to III. In each of
the three setups various algebraic structures have been
constructed in \cite{Gr, N2, S2, W} such as Hopf algebra, vertex
operators, and Heisenberg algebra, etc. We remark that there has
been a construction of an additive isomorphism between the spaces
in I and II due to de Cataldo and Migliorini \cite{CM}.

In a most important case when $X$ is $\C^2$ acted upon by a finite
subgroup $\G \subset SL_2(\C)$ and $Y$ is the minimal resolution
$\ale$ of $\C^2 /\G$, the above diagram reduces to the following
one:
\begin{eqnarray}   \label{eq_mckay}
 \begin{array}{ccc}
   \mbox{I} &  & \mbox{II}    \\
\bigoplus_{n \geq 0} H( \ale^{[n]})
          & \leftarrow -\rightarrow
          &  \bigoplus_{n \geq 0} K^{\topo}_{S_n } (\ale^n) \bigotimes \C \\
           & \; \,\nwarrow \quad \quad &\uparrow  \\
           & \quad \searrow & \downarrow  \\
          & &  \mbox{III}  \\
          & & \bigoplus_{n \geq 0} R (\Gn )
 \end{array}
\end{eqnarray}
by using the Thom isomorphism between $K^{\topo}_{\Gn} (\C^{2n})$
and the representation ring $R_{\Z}(\Gn)$ of the wreath product.
Here $R (\Gn ) = R_{\Z} (\Gn) \otimes \C$ in our notation. It was
pointed out in \cite{W} that the Frenkel-Kac-Segal homogeneous
vertex representation can be realized in terms of representation
rings of such wreath products. Such a finite group theoretic
construction, which can be viewed as a new form of McKay
correspondence, has been firmly established recently in our work
\cite{FJW} jointly with I.~Frenkel and Jing. It remains a big
puzzle however why there are many parallel algebraic structures in
these different setups.

In the present paper we propose a coherent approach to fill in the
gap (at least in the setup of diagram (\ref{eq_mckay}) above) and
present several canonical ingredients in our approach. More
explicitly, we provide direct links from wreath products to
Hilbert schemes. We find a natural interpretation of a main
ingredient (the so-called weighted bilinear form) in \cite{FJW},
which brings us one step closer to a {\em direct} isomorphism of
the two forms of McKay correspondence respectively in terms of
Hilbert schemes and wreath products. We also establish
isomorphisms of various algebraic structures in II and III. Let us
discuss in more detail.

Given a finite subgroup $\G$ of $SL_2 (\C)$, we observe that there
is a natural identification between $\C^{2n} /\Gn$ and the $n$-th
symmetric product $(\C^2 /\G)^{(n)}$ of the simple singularity
$\C^2 /\G$. The following commutative diagram \cite{W}
\begin{eqnarray*}
 \begin{array}{ccc}
   \ale^{[n]}        & \stackrel{\pi_n}{\longrightarrow}
        & (\ale)^n / S_n    \\
   \downarrow \tau_n &                 & \downarrow  \tau_{(n)}      \\
   \C^{2n} /\Gn      &   \cong            & (\C^2 / \G)^n /S_n .
 \end{array}
\end{eqnarray*}
defines a resolution of singularities $\tau_n : \alen \rightarrow
\C^{2n} /\Gn$, where $\tau_{(n)}$ is naturally induced from the
minimal resolution $\ale \rightarrow \C^2 / \G$. We show that
$\tau_n$ is a semismall crepant resolution, which provides the
first direct link between wreath products and Hilbert schemes. We
show that the fiber of $\tau_n$ over $[0]\in \C^{2n} /\Gn$
(associated to the origin of $\C^{2n}$) is of pure dimension $n$
and we give an explicit description of its irreducible components.

We conjecture that there exists a canonical isomorphism between
the Hilbert quotient $\hilqgn$ of $\C^{2n}$ by $\Gn$ (see
\cite{Ka} or Subsection~\ref{subsect_link} for the definition of
Hilbert quotient) and the Hilbert scheme $\ale^{[n]}$, and provide
several nontrivial steps toward establishing this conjecture. More
explicitly, we first single out a distinguished nonsingular
subvariety $\hnew$ of $(\C^2)^{[nN]}$ and construct a morphism
$\varphi$ from $\hilqgn$ to $\hnew$, where $N$ is the order of the
group $\G$. We use here a description of a set of generators for
the algebra of $\Gn$ invariant regular functions on $\C^{2n}$
which is a generalization of a theorem of Weyl \cite{Wey}. It
follows by construction that our morphism from $\hilqgn$ to
$\hnew$ when restricted to a certain Zariski open set is indeed an
isomorphism. We give a quiver variety description of $\hnew$ and
$\C^{2n} /\Gn$ in the sense of Nakajima \cite{N, N3}. Such an
identification follows easily from Nakajima's quiver
identification of Hilbert scheme of points on $\C^2$ (cf.
\cite{N2} and Varagnolo-Vasserot \cite{VV}). According to Nakajima
\cite{N4}, it can be shown essentially by using
Kronheimer-Nakajima \cite{KN} that the Hilbert scheme $\ale^{[n]}$
is a quiver variety associated to the same quiver data as $\hnew$
but with a different stability condition. It follows that $\hnew $
and $\alen$ are diffeomorphic by Corollary 4.2 in \cite{N}. In
this way we have obtained a second direct link between $\Gn$ and
$\ale^{[n]}$. One plausible way to establish our main conjecture
will be to establish that $\varphi$ is an isomorphism between
$\hilqgn$ and $\hnew$ and that the diffeomorphism between $\hnew $
and $\alen$ is indeed an isomorphism as complex varieties.

Our construction contains two distinguished cases which have been
studied by others. The morphism above for $n =1$ becomes an
isomorphism due to Ginzburg-Kapranov (unpublished) and
independently Ito-Nakumura \cite{INr}. Our morphism above also
generalizes Haiman's construction \cite{H} of a morphism from
$\hilqn$ to $(\C^2)^{[n]}$ which corresponds to our special case
for $\G$ trivial (where no passage to the quiver variety
description is needed). Haiman \cite{H} has in addition shown that
the morphism being an isomorphism is equivalent to the validity of
the remarkable $n!$ conjecture due to Garsia and Haiman \cite{GH}.
(We remark that there has been also attempt by Bezrukavnikov and
Ginzburg \cite{BG} in establishing this conjectural isomorphism
for $\G$ trivial.) Very recently a proof of the $n!$ conjeture
(and this isomorphism conjecture) has been announced by Haiman in
his homepage by establishing a Cohen-Macaulay property of a
certain universal scheme which was conjectured in \cite{H}. It is
natural for us to conjecture similarly a Cohen-Macaulay property
of a certain universal scheme in our setup
(Conjecture~\ref{conj_cohen}), which is sufficient to imply that
$\varphi$ is an isomorphism.

A distinguished virtual character of $\Gn$ has been used to
construct a semipositive definite symmetric bilinear form (called
a weighted bilinear form) on $R_{\Z}(\Gn)$ which plays a
fundamental role in the wreath product approach to McKay
correspondence \cite{FJW}. Indeed it is given by the $n$-th tensor
of the McKay virtual character $\la (\C^2)$ of $\G$. On the other
hand, the virtual character $\lambda (\C^{2n})$ of $\Gn$ induced
from the Koszul-Thom class defines a canonical bilinear form on
the Grothendieck group $K^0_{\Gn} (\C^{2n})$ of the bounded
derived category $D^0_{\Gn} (\C^{2n})$ consisting of $\Gn$
equivariant coherent sheaves whose cohomology sheaves are
concentrated on the origin. Although they are defined very
differently, these two virtual characters of $\Gn$ are shown to
coincide. This establishes an isometry between $K^0_{\Gn}
(\C^{2n})$ and $R_{\Z}(\Gn)$ endowed with the weighted bilinear
form, and thus provides a natural explanation of the weighted
bilinear form introduced from a purely group theoretic
consideration. (Actually we establish the coincidence of virtual
characters for more general $\G \subset GL_k (\C)$ and the induced
$\Gn$-action on $\C^{kn}$.) We regard this isometry as an
important ingredient toward a direct isomorphism of the two forms
of McKay correspondence realized respectively in terms of Hilbert
schemes \cite{N2, Gr} and of wreath products \cite{FJW}.

While our motivation is quite different, our main conjecture fits
into the scheme of Reid \cite{R} who asks for what finite subgroup
$G \subset SL_K(\C)$ the Hilbert quotient $\C^K // G$ (also called
$G$-Hilbert scheme on $\C^K$) is a crepant resolution of $\C^K/G$.
Note that the notion of McKay correspondence is meant in the
strict sense in this paper while the McKay correspondence in
\cite{R} is in a generalized sense. Our work provides supporting
evidence for an affirmative answer of the McKay correspondence in
the sense of \cite{R} for $\C^{2n}$ acted upon by $\Gn$ which in
turn is a key step to a direct isomorphism of the two form of
McKay correspondence mentioned above.

By applying a remarkable theorem of Bridgeland-King-Reid
\cite{BKR} to our situation, our main conjecture on the
isomorphism between Hilbert quotients and Hilbert schemes implies
that the equivalence of bounded derived categories among $D_{\Gn}
(\C^{2n})$, $D_{S_n} (\ale^n ),$ and $D (\ale^{[n]})$. For $n =1$
this is a theorem due to Kapranov-Vasserot \cite{KV}. Such an
equivalence can be viewed as a direct connection among the objects
in the diagram (\ref{eq_mckay}), where K-groups of topological
vector bundles are replaced by K-groups of sheaves and connection
between K-group and homology is made via Chern character.

In the end we establish a direct isomorphism of various algebraic
structures in II and III in the diagram (\ref{eq_master}). More
explicitly, we construct Schur bases of the equivariant K-groups
in II and III and show that a canonical one-to-one correspondence
between these two bases gives the desired isomorphism for various
algebraic structures such as Hopf algebras, $\lambda$-rings, and
Heisenberg algebra, etc.

The plan of the paper goes as follows. In
Sect.~\ref{sect_morphism} we study the resolution of singularities
$\tau_n$ and provide various steps toward establishing our
conjecture on the isomorphism between the Hilbert quotient
$\hilqgn$ and the Hilbert scheme $\alen$. In
Sect.~\ref{sect_mckay} we show that two virtual characters of
$\Gn$ arising from different setups coincide with each other and
discuss various implications of this and our main conjecture. In
Sect.~\ref{sect_ktheory} we establish isomorphisms of various
algebraic structures in II and III of diagram (\ref{eq_master}).

\noindent {\bf Acknowledgment.}  Some ideas of the paper were
first conceived when I was visiting the Max-Planck Institut f\"ur
Mathematik (MPI) at Bonn in 1998. I thank MPI for a stimulating
atmosphere. Hiraku Nakajima's papers and lecture notes are sources
of inspiration. I am grateful to him for the influence of his
work. His comments on an earlier version of this paper helps much
to clarify the exposition of Subsection~\ref{subsect_quiver}. I
thank Mark de Cataldo for helpful discussions and comments. I also
thank Igor Frenkel for his comments and informing me that he has
been recently thinking about some related subjects from somewhat
different viewpoint.
\section{Wreath products and Hilbert schemes}
\label{sect_morphism}
In this section, we establish some direct
connections between wreath products and Hilbert schemes, i.e.
between I and III in the diagram (\ref{eq_mckay}). In particular
for $\G$ trivial it reduces to relating I and II.
\subsection{The wreath products}
Given a finite group $\G$, we denote by $\G^*$ the set of all the
inequivalent complex irreducible characters $\{ \g_0, \g_1,
\ldots, \g_r \}$ and by $\G_*$ the set of conjugacy classes. We
denote by $\g_0$ the trivial character and by $\G^*_0$ the set of
non-trivial characters $\{\g_1, \ldots, \g_r \}$. The $\C$-span of
$\g \in \G^*$, denoted by $ R(\Gamma) $, can be identified with
the space of class functions on $\Gamma$. We denote by $\Rz$ the
integral span of irreducible characters of $ \G$.

Let $\Gamma^n = \Gamma \times \cdots \times \Gamma$ be the $n$-th
direct product of $\Gamma$. The symmetric group $S_n$ acts on
$\Gamma^n$ by permutations: $\sigma (g_1, \cdots, g_n)
  = (g_{\sigma^{ -1} (1)}, \cdots, g_{\sigma^{ -1} (n)}).
$
The wreath product of $\Gamma$ with $S_n$ is defined to be the
semi-direct product $$
 \Gamma_n = \{(g, \sigma) | g=(g_1, \cdots, g_n)\in {\Gamma}^n,
\sigma\in S_n \} $$
 with the multiplication given by
 $ (g, \sigma)\cdot (h, \tau)=(g \, {\sigma} (h), \sigma \tau ) .
$ Note that $\G^n$ is a normal subgroup of $\Gn$.

Let $\la=(\la_1, \la_2, \cdots, \la_l)$ be a partition: $\la_1\geq
\dots \geq \la_l \geq 1$. The integer $|\la|=\la_1+\cdots+\la_l$
is called the {\em weight}, and $l (\la ) =l$ is called the {\em
length} of the partition $\la $. We will often make use of another
notation for partitions: $ \la=(1^{m_1}2^{m_2}\cdots) , $ where
$m_i$ is the number of parts in $\la$ equal to $i$.

Given a family of partitions $\rho=(\rho(x))_{x\in S}$ indexed by
a finite set $S$, we define the {\em weight} of $\rho$ to be
$$\|\rho\|=\sum_{x\in S}|\rho(x)|.$$ Sometimes it is convenient to
regard $\rho=(\rho(x))_{x\in S}$ as a partition-valued function on
$S$. We denote by ${\cal P}(S)$ the set of all partitions indexed
by $S$ and by ${\cal P}_n(S)$ the set of all partitions $\rho$ in
${\cal P}(S)$ of weight $n$.

The conjugacy classes of ${\Gamma}_n$ can be described in the
following way. Let $x=(g, \sigma )\in {\Gamma}_n$, where $g=(g_1,
\cdots, g_n) \in {\Gamma}^n,$ $ \sigma \in S_n$. The permutation
$\sigma $ is written as a product of disjoint cycles. For each
such cycle $y=(i_1 i_2 \cdots i_k)$ the element $g_{i_k} g_{i_{k
-1}} \cdots g_{i_1} \in \Gamma$ is determined up to conjugacy in
$\Gamma$ by $g$ and $y$, and will be called the {\em
cycle-product} of $x$ corresponding to the cycle $y$. For any
conjugacy class $c$ and each integer $i\geq 1$, the number of
$i$-cycles in $\sigma$ whose cycle-product lies in $c$ will be
denoted by $m_i(c)$. Denote by $\rho (c)$ the partition $(1^{m_1
(c)} 2^{m_2 (c)} \ldots )$. Then each element $x=(g, \sigma)\in
{\Gamma}_n$ gives rise to a partition-valued function $( \rho
(c))_{c \in \G_*} \in {\mathcal P} ( \G_*)$ such that $\sum_{i, c}
i m_i(c) =n$. The partition-valued function $\rho =( \rho(c))_{ c
\in G_*} $ is called the {\em type} of $x$. It is well known (cf.
\cite{M,Z}) that any two elements of ${\Gamma}_n$ are conjugate in
${\Gamma}_n$ if and only if they have the same type.
\subsection{A resolution of singularities} \label{subsect_fiber}
Let $X$ be a smooth complex algebraic variety acted upon by a
finite group $\G$ of order $N$. We denote by $X^{[n]}$ the Hilbert
scheme of $n $ points on $X$ and denote by $X^{(n)} = X^n /S_n$
the $n$-th symmetric product. Both $X^{[n]}$ and $X^{(n)}$ carry
an induced $\G$-action from $X$.

Now assume $X$ is a quasi-projective surface. A beautiful theorem
of Fogarty \cite{Fo} states that $X^{[n]}$ is non-singular of
dimension $2n$. It is well known (cf. e.g. \cite{N2}) that the
Hilbert-Chow morphism $X^{[n]} \rightarrow X^{(n)}$ is a
resolution of singularities. Given a partition $\nu$ of $n$ of
length $l$: $\nu_1 \geq \nu_2 \geq \ldots \geq \nu_l
>0$, we define
\[
X^{(n)}_{\nu} = \left\{ \sum_{i =1}^l \nu_i x_i \in X^{(n)} |x_i
\neq x_i \mbox{ for } i \neq j \right\}.
\]
A natural stratification of $X^{(n)}$ is given by
\begin{eqnarray*}
X^{(n)} = \bigsqcup_{\nu} X^{(n)}_{\nu}.
\end{eqnarray*}

In the remainder of this section, we let $\G$ be a finite subgroup
of $ SL_2 (\C )$ unless otherwise specified. The classification of
finite subgroups of $SL_2 (\C)$ is well known. The following is a
complete list of them: the cyclic, binary dihedral, tetrahedral,
octahedral and icosahedral groups. We denote by $\tau : \ale \rightarrow \C^2 /\G$ the minimal resolution of the simple singularity.

A canonical identification between $\C^{2n} /\Gn$ and $(\C^2
/\G)^n /S_n $ is given as follows: given a $\Gn$-orbit, say $\Gn .
(x_1, \ldots, x_n)$ for some $(x_1, \ldots, x_n) \in (\C^2)^n =
\C^{2n}$, we obtain a point $[\G.x_1] +\ldots + [\G.x_n]$ in
$(\C^2 / \G)^n /S_n $, where $\G.x_i$ denotes the $\G$-orbit of
$x_1$, i.e. a point in $\C^2 /\G$. It is easy to see that this map
is independent of the choice of the representative
$(x_1,\ldots,x_n)$ in the $\Gn$-orbit and it is one-to-one. The
following commutative diagram \cite{W}
\begin{eqnarray}   \label{eq_mine}
 \begin{array}{ccc}
   \alen        & \stackrel{\pi_n}{\longrightarrow}
        & (\ale)^n / S_n    \\
   \downarrow \tau_n &                 & \downarrow  \tau_{(n)}      \\
   \C^{2n} /\Gn      &   \cong            & (\C^2 / \G)^n /S_n .
 \end{array}
\end{eqnarray}
defines a morphism $\tau_n : \alen \rightarrow  \C^{2n} /\Gn.$

\begin{proposition}  \label{prop_resol}
The morphism $\tau_n : \ale^{[n]} \rightarrow \C^{2n} /\Gn $ is a
semismall crepant resolution of singularities.
\end{proposition}
\begin{demo}{Proof}
It is clear by definition that $\tau_n$ is a resolution of
singularities. We now describe a stratification of $\C^{2n} /\Gn$.

The simple singularity $\C^2 /\G$ has a stratification given by
the singular point $o$ and its complement denoted by $\C_0^2 /\G$.
It follows that a stratification of $(\C^2/ \G)^n$ is given by $n
+1 $ strata $(\C^2/ \G)^n [i]$ $( 0 \leq i \leq n)$, where $(\C^2/
\G)^n [i]$ consisting of points in the Cartisian product $(\C^2/
\G)^n$ which has exactly $n - i$ components given by the singular
point. Then a stratification of $\C^{2n} /\Gn = (\C^2 / \G)^n /S_n
$ is given by
\begin{eqnarray*}
(\C^2 / \G)^n /S_n
 & =& \bigsqcup_{i=0}^n (\C^2/ \G)^n [i] /S_n \\
 & \cong & \bigsqcup_{i=0}^n (\C_0^2 /\G)^i /S_i \times \{ (n -i) o \}  \\
 & \cong & \bigsqcup_{i=0}^n (\C_0^2 /\G)^{(i)} \\
 & \cong & \bigsqcup_{i=0}^n \bigsqcup_{|\mu| =i} (\C_0^2 /\G)^{(i)}_{\mu}.
\end{eqnarray*}
The codimension of the strata $(\C_0^2/\G)^{(i)}_{\mu}$ is $2n - 2
l (\mu)$. Clearly we also have
\begin{eqnarray*}
\tau_n^{-1} ( (\C_0^2 /\G)^{(i)}_{\mu} \times \{ (n -i) o \} ) =
 (\C_0^2 /\G)^{[i]}_{\mu} \times \tau^{-1} (0)^{n -i}.
\end{eqnarray*}
It follows that the dimension of a fiber over a point in this
strata is equal to $(i - l (\mu) ) + (n -i) = n - l (\mu)$ which
is half of the codimension of the strata above. Thus $\tau_n$ is
semismall.

The canonical bundles over $\C^{2n} /\Gn$ and $\alen$ is trivial
due to the existence of holomorphic symeplectic forms (note that
$\Gn$ preserves the symeplectic form on $\C^{2n}$). Thus $\tau_n$
is crepant.
\end{demo}

\begin{remark} \rm
\begin{enumerate}
\item $\tau_n$ is a symplectic resolution in the sense that
the pullback of the holomorphic symplectic form on $\C^{2n} / \Gn$
outside of the singularities can be extended to a holomorphic
symplectic form on $\ale^{[n]}$. When $n =1$ this becomes the
minimal resolution $\tau = \tau_1 : \ale \rightarrow \C^2 /\G$.
\item  It follows from the semismallness of $\tau_n$
that $\tau_{(n)} $ is also semismall. Then the diagram
(\ref{eq_mine}) is remarkable in that all three maps $\tau_n, \pi_n$
and $\tau_{(n)} $ are semismall.
\item The resolution $\tau_n : \ale^{[n]} \rightarrow \C^{2n} /\Gn $
is one-to-one over the non-singular locus of the orbifold $\C^{2n}
/ \Gn$ corresponding to regular $\Gn$-orbits in $\C^{2n}$.
\end{enumerate}
\end{remark}

We denote by $\ale^{[n], 0}$ the fiber of $\tau_n$ over $[0]$,
where $[0]$ denotes the image in $\C^{2n}/\Gn$ of the origin of
$\C^{2n}$. We assume that $\G$ is not trivial. By the diagram
(\ref{eq_mine}), we have
\begin{eqnarray*}
\ale^{[n], 0} = \pi_n^{ -1} \tau_{(n)}^{ -1} (0) = \pi_n^{ -1}
(D^{(n)}),
\end{eqnarray*}
where $D = \tau^{-1} (0)$ is the
exceptional divisor in $\ale$. It is well known that the
irreducible components of $D$ are projective lines $\Sigma_{\g}$
parameterized by the set of non-trivial characters $\g \in \G^*_0$
(cf. e.g. \cite{GSV}).

Recall that given an irreducible curve $\Sigma \subset \ale$, the
variety
\[
L^n \Sigma := \left\{ I \in \alen | \mbox{Supp}( {\cal O} /I
)\subset \Sigma \right\} = \pi_n^{-1} (\Sigma^{(n)}).
\]
introduced by Grojnowski \cite{Gr} (also see \cite{N2}) plays an
important role in understanding the middle-dimensional homology
groups of Hilbert schemes and in connection with symmetric
functions. One can show (cf. {\em loc. cit.}) that the irreducible
components of $L^n \Sigma$ are parameterized by partitions $\nu$
of $n$ and given by
\begin{eqnarray*}
L^{\nu} \Sigma = \overline{ \pi_n^{-1} (\Sigma^{(n)}_{\nu} ) } ,
\end{eqnarray*}
where $ \Sigma^{(n)}_{\nu} $  is the stratum of the symmetric
product $\Sigma^{(n)}$ associated to $\nu$.

It is interesting to observe that the fiber $\ale^{[n], 0}$ is a
natural generalization of the above construction. Given $\rho =
(\rho (\g))_{\g \in \G^*_0} \in {\cal P}_n (\G^*_0)$, we set $n_i
= |\rho (\g_i)|$ and define
\begin{eqnarray*}
L^{\rho}D = \overline{ \pi_n^{-1}
\left((\Sigma_{\g_1})^{(n_{1})}_{\rho (\g_1)} \times \ldots \times
(\Sigma_{\g_r})^{ (n_{r})}_{\rho (\g_r)} \right)} .
\end{eqnarray*}

\begin{proposition} \label{prop_fiber}
Let $\G$ be non-trivial. The fiber $\ale^{[n], 0}$ is of pure
dimension $n$, and its irreducible components are given by
$L^{\rho}D, \rho \in {\cal P}_n (\G^*_0)$.
\end{proposition}

\begin{demo}{Proof}
The component $L^{\rho}D$ is irreducible since the fiber of
$\pi_n$ is so. $L^{\rho}D$ is of dimension $n$ since the dimension
of $(\Sigma_{\g})^{ |\rho (\g)|}_{\rho (\g)}$ $(\g \in \G^*_0)$ is
$l( \rho (\g))$ and the dimension of the fiber of $\pi_n$ is equal
to $$\sum_{i =1}^r \left( |\rho (\g_i) | - l (\rho (\g_i)) \right)
= n - \sum_i l( \rho (\g_i)). $$ Therefore the dimension of
$L^{\rho}D$ is $( n - \sum_i l( \rho (\g_i)) ) + \sum_i l( \rho
(\g_i)) = n$.
\end{demo}
\subsection{Hilbert quotient and a subvariety of $(\C^2)^{[nN]}$}
\label{subsect_link}
Let $X$ be a smooth complex algebraic variety acted upon by a
finite group $\G$ of order $N$. A regular $\G$-orbit can be viewed
as an element in the Hilbert scheme $X^{[N]}$ of $N$ points in
$X$. The Hilbert quotient is the closure $\hilq$ of the set of
regular $\G$-orbits in $X^{[N]}$ (cf. \cite{Ka}). It follows that
there exists a tautological vector bundle over $\hilq$ of rank
$N$. The group $\G$ acts on the tautological bundle fiberwise and
each fiber is isomorphic to the regular representation of $\G$.

Note that the wreath product $\Gn $ acts faithfully on the affine
space $\C^{2n} = (\C^2)^n$. We make the following conjecture which
provides a key connection between I and III in the diagram
(\ref{eq_mckay}). We remark that a more general conjecture can be
easily formulated in the setup of the diagram (\ref{eq_master}).

\begin{conjecture} \label{conj_main}
There exists a canonical isomorphism between the Hilbert quotient
$\hilqgn$ and the Hilbert scheme $\ale^{[n]}$ of $n$ points on
$\ale$.
\end{conjecture}
\begin{remark} \rm \label{rem_one}
 \begin{enumerate}
\item
Conjecture~\ref{conj_main} for $n =1$ reduces to a theorem due to
Ginzburg-Kapronov (unpublished) and independently Ito-Nakamura
\cite{INr} which says $\C^2 // \G \cong \ale$.
\item When $\G$ is trivial and so $\Gn $ is the
symmetric group $S_n$,  Haiman \cite{H} has shown that
Conjecture~\ref{conj_main} for $\G$ trivial is equivalent to a
remarkable {\em $n!$ conjecture} due to Garsia and Haiman
\cite{GH}. A proof of the $n!$ conjecture has been very recently
posted by Haiman in his UCSD webpage.
\item Conjecture~\ref{conj_main} implies an
isomorphism of Hilbert quotients:
\[
 \hilqgn \cong (\C^2 // \G )^n //S_n,
\]
which seems to be in a more symmetric form. On the other hand,
assuming the $n!$-conjecture, one can show that an isomorphism of
Hilbert quotients above implies Conjecture~\ref{conj_main}.
\end{enumerate}
\end{remark}

Since $\tau_n : \alen \rightarrow \C^{2n } /\Gn$ is a resolution
of singularities ( and thus is proper) and $\C^{2n } /\Gn$ is a
normal variety, the algebra of regular functions on $\alen$ is
isomorphic to the algebra of $\Gn$ invariants on the regular
functions on $\C^{2n}$. The following lemma gives a description of
the algebra of $\Gn$ invariants in $\C [\x, \y]$, where we denote
by $\x$ (resp. $\y$) the $n$-tuple $x_1, \ldots, x_n$ (resp. $y_1,
\ldots, y_n$). It generalizes Weyl's theorem \cite{Wey} for
symmetric groups and the proof is similar.

\begin{lemma} \label{lem_weyl}
The algebra of invariants $\C [\x, \y]^{\Gn}$ is generated by
\begin{eqnarray*}
\tilde{f} (\x, \y) = f(x_1, y_1) + f(x_2, y_2)
                     + \ldots + f(x_n, y_n),
\end{eqnarray*}
where $f$ runs over an arbitrary linear basis $\cal B$ for
the space of invariants $\C [x,y]^{\G}$.
\end{lemma}

\begin{demo}{Proof}
We prove the lemma by induction on $n$. When $n =1$ it is evident.
Assume now that we have established the lemma for $n-1$. We use
$\x '$ and $ \y'$ to denote $x_2, \ldots, x_n$ and respectively
$y_2, \ldots, y_n$. The space $\C[ \x ', \y ' ]$ is acted on by
the wreath product subgroup $\G_{n -1} \subset \Gn$.

Given any $F (\x, \y) \in \C [\x, \y]^{\Gn}$, we can write it as a
linear combination of $x_1^{\alpha} y_1^{\beta} F_{\alpha
\beta}'(\x ', \y '),$ where $ F_{\alpha \beta}'(\x ', \y ')$ is
some $\G_{n -1}$-invariant polynomial. By induction assumption, we
can write $ F_{\alpha \beta}'(\x ', \y ')$ as a polynomial in
terms of $\tilde{f} (\x', \y')$ where $f \in {\cal B}$, and in
turn as a polynomial in terms of $\tilde{f} (\x, \y)$ and $x_1,
y_1$. Therefore we can write $F (\x, \y)$ as a linear combination
of polynomials of the form $ x_1^{\alpha} y_1^{\beta} G_{\alpha
\beta}(\x, \y),$ where $ G_{\alpha \beta}(\x, \y )$ is a
polynomial in terms of $\tilde{f} (\x', \y')$ where $f \in \cal
B$. Since both $ F_{\alpha \beta}(\x, \y)$ and $G_{\alpha
\beta}(\x, \y)$ are $\Gn$-invariant and thus in particular
invariant with respect to the symmetric group $S_n$ and the first
factor $\G$ in $\G^n \subset \Gn$, $F(\x, \y)$ becomes a linear
combination of $p_{\alpha \beta}(\x, \y) G_{\alpha \beta}(\x,\y).$
Here $p_{\alpha \beta}(\x, \y)$ denotes $\frac1{nN} \sum_{i =1}^n
\sum_{g \in \G} (g. x_i)^{\alpha} (g.y_i)^{\beta}$, the average of
$x_1^{\alpha} y_1^{\beta}$ over $\G \times S_n$ (which is the same
as the average over $\Gn$). We complete the proof by noting that
$p_{\alpha \beta}(\x, \y)$ is a linear combination of $\tilde{f}
(\x, \y)$ where $f \in \cal B$.
\end{demo}

Given $J \in \hilqgn$ we regard it as an ideal in $\C [\x, \y]$ of
colength $N^n n!$ (which is the order of $\Gn$). Then the quotient
$\C [\x, \y ] / J$ affords the regular representation $R$ of
$\Gn$, and its only $\Gn$-invariants are constants. Thus we have
$\tilde{f} (\x, \y ) = c_f \mbox{ mod }J$ for some constant $c_f$.
Recall that $\G_{n -1}$ acts on $\C [\x ', \y ']$. By
Lemma~\ref{lem_weyl}, the space $\C [\x, \y]^{\G_{n -1}}$ is
generated by $x_1, y_1$ and $f(x_2, y_2) + \ldots + f(x_n, y_n)$.
The latter is equal to $c_f - f(x_1, y_1) \mbox{ mod } J$. Thus
$(\C [\x, \y] / J)^{\G_{n -1}}$ is generated by $x_1, y_1$ and
$c_f - f(x_1, y_1), f \in \C [x, y]^{\G}$. It follows that
\begin{eqnarray}  \label{eq_link}
  \C [x_1, y_1 ] / (J \cap \C [x_1, y_1 ])
  \equiv ( \C [\x, \y] / J)^{\G_{n -1}},
\end{eqnarray}
which has dimension $nN = | \Gn | / | \G_{n -1} |$ because
$(\C[\x, \y] /J)^{\G_{n -1}}$ can be identified with the space of
${\G_{n -1}}$-invariants in the regular representation of $\G$.

The first copy of $\G$ in the Cartesian product $\G^n \subset \Gn$
commutes with $\G_{n -1}$ above. It follows from (\ref{eq_link})
that the quotient $\C [x_1, y_1 ] / (J \cap \C [x_1, y_1 ])$ as a
$\G$-module is isomorphic to $R^n$, a direct sum of $n$ copies of
the regular representation $R$ of $\Gn$.

The $\G$-action on $\C^2$ induces a $\G$-action on the Hilbert
scheme $(\C^2 )^{[nN]}$ and the symmetric product $
(\C^2)^{(nN)}$. The Hilbert-Chow morphism $(\C^2 )^{[nN]}
\rightarrow (\C^2 )^{(nN)}$ induces one between the set of
$\G$-fixed points $(\C^2 )^{[nN], \G} \rightarrow (\C^2
)^{(nN),\G}$. As the fixed point set of a non-singular variety by
the action of a finite group, $(\C^2 )^{[nN], \G}$ is
non-singular. Denote by $\hnew$ the set of $\G$-invariant ideals
$I$ in the Hilbert scheme $ (\C^2)^{[nN]}$ such that the quotient
$\C [x, y] / I$ is isomorphic to $R^n$ as a $\G$-module. Since the
quotient $\C [x, y]/I$ are isomorphic as $\G$-modules for all $I$
in a given connected component of $(\C^2 )^{[nN], \G}$, the
variety $\hnew$ is a union of components of $(\C^2 )^{[nN], \G}$.
In particular $\hnew$ is non-singular (we shall see that $\hnew$
is indeed connected of dimension $2n$).

Therefore by sending the ideal $J$ to the ideal $J \cap \C [x_1,
y_1]$, we have defined a map $\varphi$ from $ \hilqgn$ to $ (\C^2
)^{[nN]}$, whose image lies in $\hnew$. We also denote $\varphi:
\hilqgn \longrightarrow \hnew$.

The map $\varphi$ can be also understood as follows. Let
$\quotuniv$ be the universal family over the Hilbert quotient
$\hilqgn$ which is a subvariety of the Hilbert scheme
$(\C^{2n})^{[n!N^n]}$:
\begin{eqnarray}  \label{eq_quotuniv}
 \begin{array}{ccc}
  \quotuniv & \longrightarrow & \C^{2n}    \\
   \downarrow   &                 &        \\
    \hilqgn &        &
 \end{array}
\end{eqnarray}
It has a natural $\Gn$-action fiberwise such that each fiber
carries the regular representation of $\Gn$. Then $\quotuniv/
\G_{n-1}$ is flat and finite of degree $nN$ over $\hilqgn$, and
thus can be identified with a family of subschemes of $\C^2$ as
above. Then $\varphi$ is the morphism given by the universal
property of the Hilbert scheme $(\C^2 )^{[nN]}$ for the family
$\quotuniv/ \G_{n-1}$.

Thus we have established the following.

\begin{theorem} \label{th_morph}
We have a natural morphism $\varphi : \hilqgn \longrightarrow
\hnew$ defined as above.
\end{theorem}

\begin{remark} \rm
For $\G $ trivial, $\Gn$ reduces to $S_n$ and $\hnew$ becomes the
Hilbert scheme $(\C^2)^{[n]}$. In this case the above morphism
$\hilqn \rightarrow (\C^2)^{[n]}$ was earlier constructed by
Haiman \cite{H}. We plan to elaborate further on generalizations
of \cite{H} to our setup in the future.
\end{remark}

\begin{conjecture}
The morphism $\varphi : \hilqgn \longrightarrow \hnew$ is an
isomorphism.
\end{conjecture}

Observe that $ (\C^2 )^{(nN),\G}$ is the subset of $(\C^2)^{(nN)}$
consisting of points
\begin{eqnarray*}
\sum_{\g \in \G} [\g \cdot x_1] + \ldots + \sum_{\g \in \G} [\g \cdot
x_a] + (n -a) N [0], \quad 1 \leq a \leq n, x_1 , \ldots , x_a \in
\C^2 \backslash 0,
\end{eqnarray*}
which can be thought as the $\Gn$-orbit of $(x_1, \ldots, x_a, 0,
\ldots, 0) \in (\C^2)^n = \C^{2n}$. In this way $(\C^2)^{(nN),\G}$
is identified with $\C^{2n} / \Gn$. Thus we have proved the
following proposition  by noting the inclusion $\hnew \subset
(\C^2 )^{[nN], \G} $.

\begin{proposition}  \label{prop_fix}
The $\G$-fixed point set $(\C^2)^{(nN),\G}$ of the symmetric
product $(\C^2)^{(nN)}$ can be canonically identified with
$\C^{2n} / \Gn$. The Hilbert-Chow morphism $(\C^2)^{[nN]}
\rightarrow (\C^2)^{(nN)}$, when restricted to the $\G$-fixed
point set, induces a canonical morphism $\hnew \rightarrow \C^{2n}
/ \Gn$.
\end{proposition}

\begin{remark}  \rm  \label{rem_gener}
Take an unordered $n$-tuple $T$ of distinct $\G$-orbits in $\C^2
\backslash 0$. Such an $n$-tuple defines a set of $nN$ distinct
points in $\C^2$, and thus can be regarded as an ideal $I(T)$ in
the Hilbert scheme $(\C^2)^{[nN]}$. This ideal is clearly
$\G$-invariant and as a $\G$-module $\C [x, y] / I(T)$ is
isomorphic to $R^n$. On the other hand observe that such an
$n$-tuple $T$ can be canonically identified with a regular
$\Gn$-orbit in $\C^{2n}$. In this way the sets of generic points
in $\hnew$ and $\C^{2n} /\Gn$ coincide. It is easy to see that the
morphism $\hnew \rightarrow \C^{2n} / \Gn$ above is surjective and
it is one-to-one over the set of generic points in $\C^{2n} / \Gn$
consisting of regular $\Gn$-orbits.
\end{remark}

We define the reduced universal scheme $W_{\G, n}$ as the reduced
fibered product
\begin{eqnarray*}
\begin{array}{ccc}
    W_{\G, n}  & {\longrightarrow} & \C^{2n} \\
  \downarrow   &                   & \downarrow \\
 \hnew     & \stackrel{ \tau_n}{\longrightarrow} & \C^{2n} /\Gn .
\end{array}
\end{eqnarray*}
It is known that a finite surjective morphism from $Z$ to a
non-singular variety is flat if and only if $Z$ is Cohen-Macaulay.

\begin{conjecture}  \label{conj_cohen}
$W_{\G, n}$ is Cohen-Macaulay.
\end{conjecture}

Under the assumption of the validity of
Conjecture~\ref{conj_cohen}, the universal properties of the
Hilbert scheme $(\C^{2n})^{[n!N^n]}$ induces a morphism $\psi:
\hnew \rightarrow (\C^{2n})^{[n!N^n]}$, whose image lies in
$\hilqgn$. By Remark~\ref{rem_gener} and the fact that the set of
generic points of $\hilqgn$ and $\C^{2n} / \Gn$ coincide, the two
morphisms $\varphi$ and $\psi$ are mutually inverse to each other
over generic points. Then it follows that they are inverse
everywhere, establishing Conjecture~\ref{conj_main}.

\begin{remark} \rm
The above conjecture~\ref{conj_cohen} for $\G$ trivial was first
conjectured by Haiman \cite{H}. The proof announced very recently
by Haiman in his UCSD homepage of $n!$ conjecture is based on a
proof of this conjecture of his.
\end{remark}
\subsection{A quiver variety description}  \label{subsect_quiver}
We first recall (cf. \cite{N2}) that the Hilbert scheme
$(\C^2)^{[K]}$ of $K$ points in $\C^2$ admits a description in
terms of a quiver consisting of one vertex and one arrow starting
from the vertex and returning to the same vertex itself. More
explicitly, we denote
\begin{eqnarray*}
 \widetilde{H}(K) =
\left\{
\begin{array}{lcl}
              &|& i) [B_1, B_2 ] +ij =0  \\
(B_1, B_2, i, j) &|& ii) \mbox{ there exists no proper subspace
}\\
              &|&  S \subset
              \C^K \mbox{ such that } B_{\alpha} (S) \subset S \mbox{ and } \\
               &|& \mbox{ im } i \subset S \quad (\alpha =1,2)
\end{array}
              \right\}  ,
\end{eqnarray*}
where $B_1, B_2 \in \End (\C^K ), i \in \Hom (\C, \C^K ), j \in
\Hom (\C^K, \C )$. Then we have an isomorphism
\begin{eqnarray}  \label{eq_quot}
(\C^2)^{[K]} \cong \widetilde{H}(K)/ GL_K (\C ) ,
\end{eqnarray}
where the action of $GL_K (\C )$ on $\widetilde{H}(K)$ is given by
\begin{eqnarray*}
  g \cdot (B_1, B_2, i, j) = (g B_1 g^{ -1}, g B_2 g^{ -1}, g i, jg^{-1}).
\end{eqnarray*}
It is also often convenient to regard $(B_1, B_2)$ to be in $\Hom
(\C^K , \C^2 \otimes \C^K )$. We remark that one may drop $j$ in
the above formulation because one can show by using the stability
condition that $j =0$ (cf. \cite{N2}).

The bijection in (\ref{eq_quot}) is given as follows. For $I \in
(\C^2)^{[K]}$, i.e. an ideal in $\C [x, y]$ of colength $K$, the
multiplication by $x, y$ induces endomorphisms $B_1, B_2$ on the
$K$-dimensional quotient $\C [x, y] /I$, and the homomorphism $i
\in \Hom (\C, \C^K)$ is given by letting $i (1) =1 \mbox{ mod }
I$. Conversely, given $(B_1, B_2, i)$, we define a homomorphism
$\C [x, y] \rightarrow \C^K$ by $f \mapsto f(B_1, B_2) i (1)$. One
can show by the stability condition that the kernel $I$ of this
homomorphism is an ideal of $\C [x, y]$ of colength $K$. One
easily checks that the two maps are inverse to each other.

Set $K =nN$, where $N$ is the order of $\G$. We may identify
$\C^K$ with $R^n$, $\C^2$ with the defining representation $Q$ of
$\G$ by the embedding $\G \subset SL_2 (\C)$, and $\C$ with the
trivial representation of $\G$. Denote by
\begin{eqnarray*}
M (n) & =& \Hom (R^n , Q \otimes R^n ) \bigoplus \Hom (\C,
R^n)\bigoplus \Hom (R^n, \C).
\end{eqnarray*}
By definition $\widetilde{H} (nN) \subset M(n)$. Let $GL_{\G} (R)$
be the group of $\G$-equivariant automorphisms of $R$. Then the
group $G \equiv GL_{\G} (R^n) \cong GL_n (\C) \times GL_{\G} (R)$
acts on the $\G$-invariant subspace $M(n)^{\G} $. We have the
following description of $\hnew$ as a quiver variety. This result
is known to Nakajima \cite{N4} (also cf. Theorem 1 of
Varagnolo-Vasserot \cite{VV})\footnote{I.~Frenkel informed us that
he also noticed this recently.}.
\begin{theorem}  \label{th_quiv}
The variety $\hnew $ admits the following description:
\begin{eqnarray*}
  \hnew \cong (\widetilde{H} ( nN) \cap M(n)^{\G} ) / GL_{\G} (R^n).
\end{eqnarray*}
It particular $\hnew $ is non-singular of pure dimension $2n$.
\end{theorem}

\begin{remark}   \rm \label{rem_clear}
Consider the $\G$-module decomposition $Q \otimes V_{\g_i}=
\bigoplus_j a_{ij} V_{\g_j}$, where $a_{ij} \in \Z_+$, and
$V_{\g_i} $ $(i =0, \ldots, r)$ are irreducible representations
corresponding to the characters $\g_i$ of $\G$. Set $\dim V_{\g_i}
= n_i$. Then
\begin{eqnarray}
 M(n)^{\G}
 & =& \Hom_{\G} (R^n ,Q \otimes R^n ) \bigoplus \Hom_{\G} (\C, R^n)
  \bigoplus \Hom_{\G} (R^n, \C)
 \label{eq_invt} \\
 & =& \Hom_{\G} (\sum_i \C^{n n_i} \otimes V_{\g_i},
  \C^2 \otimes \sum_i \C^{n n_i} \otimes V_{\g_i} )\nonumber \\
  && \bigoplus \Hom_{\G} (\C, R^n)  \bigoplus \Hom_{\G} (R^n, \C) \nonumber \\
 & =& \sum_{i j} a_{ij} \Hom ( \C^{n n_i} ,\C^{n n_j})
  \bigoplus \Hom (\C, V_{\g_0}^n ) \bigoplus \Hom (V_{\g_0}^n, \C).  \nonumber
\end{eqnarray}
where $\Hom_{\G}$ stands for the $\G$-equivariant homomorphisms.
In the language of quiver varieties as formulated by Nakajima
\cite{N, N3}, the above desciption of $\hnew$ identifies $\hnew$
with a quiver variety associated to the following data: the graph
consists of the same vertices and edges as the McKay quiver which
is an affine Dynkin diagram associated to a finite subgroup $\G$
of $SL_2 (\C)$; the vector space $V_i$ associated to the vertex
$i$ is isomorphic to the direct sum of $n$ copies of the $i$-th
irreducible representation $V_{\g_i}$; the vector space $W_i =0$
for nonzero $i$ and $W_0 =\C$.
\end{remark}

\begin{demo}{Proof of Theorem~\ref{th_quiv}}
Our proof is modeled on the proof of Theorem 1.9 and Theorem 4.4
in \cite{N2} which are special cases of our isomorphism for $\G$
trivial and for $n =1$ respectively. We sketch below for the
convenience of the reader.

One shows $j =0$ by using the stability condition. The isomorphism
statement follows directly from the description of ${\C^2}^{[nN]}$
given by (\ref{eq_quot}), the definition of $\hnew$, and
Eq.~(\ref{eq_invt}). We have seen earlier that $\hnew$ is
nonsingular by construction.

One shows by a direct check that $[B_1, B_2]$ is $\G$-invariant
endomorphism in $R^n$ for $(B_1, B_2) \in \Hom_{\G} (R^n , \C^2
\otimes R^n )$. The cokernel of the differential of the map
$(B_1,B_2, i) \mapsto [B_1, B_2]$ from $M(n)^{\G}$ to
$\End_{\G}(R^n)$ consists of the $\G$-endomorphisms in $R^n$ which
commute with $B_1 $ and $B_2$. By sending $f \mapsto f( i (1))$ we
define a map from the cokernel to the $n$-dimensional space
$\Hom_{\G} (\C, R^n) \cong V_{\g_0}^n$. Conversely, given $v \in
V_{\g_0}^n$, an endomorphism $f$ in $R^n$ is uniquely determined
by the equation $f(B_1^a B_2^bi(1)) = B_1^a B_2^b v$ by the
stability condition ii) in the definition of $\widetilde{H} (nN)$.
One further checks that $f$ lies in the cokernel. These two maps
are inverse to each other. Thus the cokernel has constant
dimension $n$.

The dimension of $\widetilde{H} ( nN) \cap M(n)^{\G}$ is equal to
$\dim M^{\G}(n) +n - \dim GL_{\G} (R^n) $ since the dimension of
the cokernel is $n$. Thus the quotient description of $\hnew$
implies that the dimension of $\hnew $ near $(B_1, B_2, i)$ is
given by
\begin{eqnarray*}
  & & \dim M(n)^{\G} +n - 2 \dim GL_{\G} (R^n) \\
  & =& (n^2 \dim \Hom_{\G} (R, \C^2 \otimes R) +n)+ n -2 n^2 \dim GL_{\G}
  (R) \\
  & =& n^2 (2 N) + n + n - 2 n^2 N \\
  & =& 2n,
\end{eqnarray*}
which is independent of which component a point $(B_1, B_2, i)$
is in. Here we
have used the fact that the (complex) dimension of $ \Hom_{\G} (R,
\C^2 \otimes R)$ is equal to $N$ (cf. Kronheimer \cite{Kr}).
\end{demo}

\begin{remark}  \rm
Recall that the minimal resolution $\ale$ endowed with certain
hyper-Kahler structures are called an ALE spaces \cite{Kr}.
According to Nakajima \cite{N4}, one can show that the Hilbert
scheme $\ale^{[n]}$ over an ALE space admits a quiver variety
description in terms of the same quiver data as specified in
Remark~\ref{rem_clear} but with a different stability condition,
by a modification of the proof for the description of the moduli
space of vector bundles over an ALE space \cite{KN}. It follows by
Corollary 4.2 of \cite{N} that $\alen$ and $\hnew$ is
diffeomorphic. We conjecture that they are indeed isomorphic as
complex varieties. In this way we would have obtained a morphism
$\varphi: \hilqgn \rightarrow \alen$ by combining with
Theorem~\ref{th_morph}.
\end{remark}

\begin{remark} \rm
It follows readily from the affine algebro-geometric quotient
description of the symmetric product $(\C^2)^{(nN)}$
(Proposition~2.10, \cite{N2}) that the $\G$-fixed-point set
$(\C^2)^{(nN),\G}$ or rather the orbifold $\C^{2n} / \Gn$ (see
Proposition~\ref{prop_fix}) has the following description (also
cf. \cite{VV}, Theorem~1):
\[
\C^{2n} / \Gn \cong \{(B_1, B_2, i, j) \in M(n)^{\G} | [B_1, B_2]
+ ij =0 \}
//GL_{\G} (R^n).
\]
It follows from the general theory of quiver varieties that there
is a natural projective morphism $\alen \rightarrow \C^{2n} / \Gn$
which is a semismall resolution. We expect that this is the same
as the semismall resolution $\tau_n :\alen \rightarrow
\C^{2n}/\Gn$ explicitly constructed by the diagram
(\ref{eq_mine}). We also expect that the intermediate variety
$(\ale)^n/S_n$ (see (\ref{eq_mine})) can also be identified with a
quiver variety associated with the same quiver data (as specified
in Remark~\ref{rem_clear}) but with a new stability condition. In
this case it follows from Corollary 4.2 in \cite{N} that the fiber
$\ale^{[n], 0}$ is a lagrangian subvariety in $\alen$ and it is
homotopy equivalent to $\alen$ (compare with
Proposition~\ref{prop_fiber}).
\end{remark}

\begin{remark} \rm
Quiver varieties are connected \cite{N, N3}. Thus the variety
$\hnew$ can also be defined as the closure of the set of ideals
$I(T)$ in $(\C^2)^{[nN]}$ associated to unordered $n$-tuples $T$
of distinct $\G$-orbits in $\C^2 \backslash 0$.
\end{remark}
\subsection{Canonical vector bundles}
Since $\ale$ is isomorphic to the Hilbert quotient $\hilqg$, there
exists the tautological vector bundle $\cal R$ on $\ale$ of rank
$N$, whose fiber affords the regular representation of $\G$ (cf.
\cite{GSV, N2}). It decomposes as follows:
\begin{eqnarray*}
{\cal V} \cong \sum_{\g \in \G^*} {\cal R}_{\g} \bigotimes V_{\g},
\end{eqnarray*}
where $V_{\g}$ is the irreducible representation of $\G$
associated to $\g$ and ${\cal R}_{\g}$ is a vector bundle over
$\ale$ of rank equal to $ \deg \g$ (by definition $\deg \g = \dim
V_{\g}$).

One can associate a vector bundle $E^{[n]}$ of rank $dn$ on the
Hilbert scheme $\alen$ to a rank $d$ vector bundle $E$ over $\ale$
as follows: let $U \subset \alen \times \ale$ be the universal
family
\begin{eqnarray*}
\begin{array}{ccc}
     U  &  \stackrel{ p_1}{\longrightarrow} & \ale \\
 p_2 \downarrow\quad  &                  & \\
 \alen &  &
\end{array}
\end{eqnarray*}
$U$ is flat and finite of degree $n$ over $\alen$. Then $E^{[n]}$
is defined to be $(p_2)_* p_1^* E$. In this way we obtain
canonical vector bundles ${\cal R}^{[n]}$ and ${\cal
R}^{[n]}_{\g}$ $(\g \in \G^*)$ over $\alen$ associated to $\cal R$
and ${\cal R}_{\g}$ above.

There exists a tautological vector bundle ${\cal R}^{\{n \} }$ of
rank $nN$ over $\hnew$ (and thus over $\ale^{[n]}$) induced from
the inclusion $\hnew \subset (\C^2)^{[nN]}$. The group $\G$ acts
on ${\cal R}^{\{n \} }$ fiberwise such that each fiber as a
$\G$-module is isomorphic to $R^n$. Then we have a decomposition:
\begin{eqnarray*}
{\cal R}^{\{n \} }
 = \bigoplus_{\g \in \G^*} {\cal R}^{\{n \} }_{\g} \bigotimes V_{\g},
\end{eqnarray*}
where ${\cal R}^{\{n \} }_{\g}$ is a vector bundle over
$\ale^{[n]}$ of rank equal to $n \deg \g$.

The principle bundle
\begin{eqnarray*}
\widetilde{H} (K) \cap M^{\G}(n) \stackrel{/ G}{\longrightarrow}
\hnew \end{eqnarray*}
(which follows from Theorem~\ref{th_quiv}) gives rise to various
canonical vector bundles associated to canonical representations
of $G \equiv GL_{\G} (R^n) \cong GL_n (\C) \times GL_{\G} (R)$.
For example, the vector bundle associated to the representation
$GL_{\G} (R^n) \hookrightarrow GL (R^n)$ is exactly the above
tautological vector bundle ${\cal R}^{\{n \} }$ over $\hnew$ (or
rather over $\alen$); the one associated to $GL_{\G} (R^n)
\rightarrow G L ( \g^n)$ is $ {\cal R}^{\{n \} }_{\g}$.

In the remainder of this section we assume the validity of
Conjecture~\ref{conj_main} that $ \hilqgn $ is isomoprhic to
$\alen$. An immediate corollary is the existence of a tautological
bundle $\cal V$ on $\alen$ whose fiber affords the regular
representation of $\Gn$. This comes from the tautological bundle
over the Hilbert quotient $\hilqgn$. It is well known (cf. e.g.
\cite{M, Z}) that the irreducible representations $S_{\rho}  $ of
$\Gn$ are parameterized by the set ${\cal P}_n (\G^*)$ of
partition-valued functions $\rho$ of weight $n$ on the set $\G^*$
of irreducible characters of $\G$. It follows that one has a
decomposition
\begin{eqnarray} \label{eq_wreathreg}
 {\cal V} = \bigoplus_{\rho \in {\cal
P}_n (\G^*)} S_{\rho} \bigotimes {\cal V}_{\rho},
\end{eqnarray}
where ${\cal V}_{\rho}$ is a vector bundle on $\ale^{[n]}$ of rank
equal to $\dim S_{\rho} $. The vector bundles ${\cal V}_{\rho}$
are expected to be a basis for the K-group of $\ale^{[n]}$.

We denote by ${\cal R}^{\langle n \rangle }$ the subbundle of
$\cal V$ which is given by the (fiberwise) $\G_{n -1}$ invariants
of $\cal V$. Since the first copy $\G$ in $\G^n \subset \Gn$
commutes with $\G_{n -1}$, ${\cal R}^{\langle n \rangle }$ has a
$\G$ action fiberwise such that each fiber affords the $\G$-module
$R^n$ (cf. Subsect.~\ref{subsect_link}). It decomposes as
\begin{eqnarray*}
{\cal R}^{\langle n \rangle } = \bigoplus_{\g \in \G^*} {\cal
R}^{\langle n \rangle }_{\g} \bigotimes V_{\g}.
\end{eqnarray*}
where $ {\cal R}^{\langle n \rangle }_{\g}$ is a vector bundle
over $\alen$ of rank $n \deg \g$.

We define the reduced universal scheme $\univ$ as the reduced
fibered product
\begin{eqnarray}  \label{eq_univ}
\begin{array}{ccc}
     \univ  & {\longrightarrow} & \C^{2n} \\
  \downarrow &                  & \downarrow \\
 \alen & \stackrel{ \tau_n}{\longrightarrow} & \C^{2n} /\Gn .
\end{array}
\end{eqnarray}
Under the isomorphism $\varphi: \hilqgn \rightarrow \alen$, the
universal schemes $\quotuniv, \univ$ defined respectively by
(\ref{eq_quotuniv}) and (\ref{eq_univ}) can be identified. The
following proposition follows now from the way we define $\varphi$
(cf. Subsect.~\ref{subsect_link}).

\begin{proposition}
There is a natural identification between the vector bundles
${\cal R}^{[n]}$ and ${\cal R}^{\langle n \rangle }$, respectively
${\cal R}^{[n]}_{\g}$ and ${\cal R}^{\langle n \rangle }_{\g}$.
\end{proposition}
\section{On the equivalence of two forms of McKay correspondence}
\label{sect_mckay}
\subsection{A weighted bilinear form}
In this subsection we recall the notion of a {\em weighted}
bilinear form on $R(\Gn)$ introduced in \cite{FJW}.

The standard bilinear form on $R(\G )$ is defined as follows:

\begin{eqnarray*}
\langle f, g \rangle_{\G} = \frac1{ | \G |}\sum_{x \in \Gamma}
          f(x) g(x^{ -1}).
\end{eqnarray*}
We will write $\langle \ , \ \rangle$ for $\langle \ , \
\rangle_{\G }$ when no ambiguity may arise.

Let us fix a virtual character $ \wt \in R(\G)$. The
multiplication in $ R(\G)$ corresponding to the tensor product of
two representations will be denoted by $* $. Recall that $\g_0,
\g_1, \ldots, \g_r$ are all the inequivalent irreducible
characters of $\G$ and $\g_0$ denotes the trivial character. We
denote by $a_{ij} \in \Z$ the (virtual) multiplicities of $\g_j$
in $ \wt * \g_i $, i.e. $ \wt  *  \g_i = \sum_{j =0}^r a_{ij}
\g_j.$ We denote by $A$ the $ (r +1) \times (r +1)$ matrix $ (
a_{ij})_{0 \leq i,j \leq r}$.

We introduce the following {\em weighted bilinear form} on
$R(\G)$:
 $$
  \langle f, g \rangle_{\wt } = \langle \wt * f ,  g \rangle_{\G },
   \quad f, g \in R( \G).
$$
It follows that $\langle \g_i, \g_j \rangle_{\wt } = a_{ij}.$

Throughout this paper we will always assume that $\wt $ is a {\em
self-dual}, i.e. $\wt (x) = \wt (x^{-1}), x \in \G$. The
self-duality of $\wt$ implies that $ a_{ij} = a_{ji}, $ i.e. $A$
is a symmetric matrix.

Given a representation $V$ of $\G$ with character $\g \in R(\G)$,
the $n$-th outer tensor product $V^{ \otimes  n} $ of $V$ can be
regarded naturally as a representation of the wreath product $\Gn$
whose character will be denoted by $\eta_n ( \g )$: the direct
product $\G^n$ acts on $\g^{\otimes n}$ factor by factor while
$S_n$ by permuting the $n$ factors. Denote by $\varepsilon_n$ the
(1-dimensional) sign representation of $\Gn$ on which $\G^n$ acts
trivially while $S_n$ acts as sign representation. We denote by
$\varepsilon_n ( \g ) \in R(\Gn)$ the character of the tensor
product of $\varepsilon_n$ and $V^{\otimes n}$.

We may extend naturally $\eta_n$ to a map from $R(\G)$ to
$R(\Gn)$. In particular, if $\beta$ and $\g $ are characters of
$\G$, then
\begin{eqnarray}  \label{eq_virt}
  \eta_n (\beta - \g) =
  \sum_{m =0}^n ( -1)^m \mbox{Ind}_{\G_{n -m} \times \G_m }^{\Gn}
   [ \eta_{n -m} (\beta) \otimes \varepsilon_m (\g ) ] .
\end{eqnarray}

We define a {\em weighted bilinear form} on $R( \Gn)$ by
letting $$
  \langle  f, g\rangle_{\wt, \Gn } =
   \langle \eta_n (\wt ) * f, g \rangle_{\Gn} ,
   \quad f, g \in R( \Gn).
$$ One can show that the bilinear form $\langle \ , \
\rangle_{\wt, \Gn}$ is symmetric.
A symmetric bilinear form on $\RG = \bigoplus_{n} R(\Gn)$ is
then given by
\[
\langle u, v \rangle_{\wt}
 = \sum_{ n \geq 0} \langle u_n, v_n \rangle_{\wt, \Gn } ,
\]
where $u = \sum_n u_n$ and $v = \sum_n v_n$ with $u_n, v_n \in
\Gn$.

We further specialize to the case when $\G$ is a finite subgroup
of $SL_2 (\C)$ by an embedding $\pi$, and fix the virtual
character $\wt$ of $\G$ to be $$
 \la (\pi) \equiv \sum_{i=0}^2 (-1)^i \Lambda^i \pi =2\g_0 - \pi.
$$ where $\Lambda^i$ denotes the $i$-th exterior power. We
construct a diagram with vertices corresponding to elements $\g_i$
in $\G^*$ and we draw one edge (resp. two edges) between the
$i$-th and $ j$-th vertices if $a_{ij} = -1$ (resp. $-2$).
According to McKay \cite{Mc}, the associated diagram can be
identified with affine Dynkin diagram of ADE type and the matrix
$A$ is the corresponding affine Cartan matrix. It is shown in
\cite{FJW} that the weighted bilinear form on $R_{\Z}(\Gn)$ is
semipositive definite symmetric.
\subsection{Identification of two virtual characters}
In this subsection we set $\G$ to be an arbitrary (not necessarily
finite) subgroup of $GL_k(\C)$ unless otherwise specified. We
denote by $\pi$ the $k$-dimensional defining representation of
$\G$ for the embedding.

The wreath product $\Gn$ acts naturally on $\C^{kn} =(\C^k)^n$ by
letting $\G^n$ act factor-wise and $S_n$ act as permutations of
$n$-factors. We denote by $\lambda (\C^{kn})$ the virtual
character $\sum_{i =0}^{kn} ( -1)^i \Lambda^i \C^{kn}$ of $\Gn$,
where the $i$-th exterior power $ \Lambda^i \C^{kn}$ carries an
induced $\Gn$-action. The geometric significance of $\lambda
(\C^{kn})$ will become clear later. Let $\eta_n (\lambda (\pi))$
be the virtual character of $\Gn$ built on the $n$-th tensor of
the virtual character $\lambda (\pi ) = \sum_{i =0}^k (-1)^i
\Lambda^i \C^k$ of $\G$.
\begin{theorem} \label{th_key}
 The virtual characters $\lambda (\C^{kn})$ and
 $\eta_n (\lambda (\pi))$ of $\Gn$ are equal.
\end{theorem}

\begin{demo}{Proof}
 Given $x \times y \in \Gn$, where
 $x \in \G_m$ and $y \in \G_{n -m}$ for some $m$, we have by definition
\begin{eqnarray}   \label{eq_mult}
 \eta_n (\lambda (\pi) ) (x \times y)
   & =& \eta_m (\lambda (\pi)) (x)\; \eta_{n-m} (\lambda (\pi))(y),
            \nonumber   \\
  \lambda (\C^{kn}) (x \times y)
   & =& \lambda (\C^{km}) (x)\; \lambda (\C^{k(n -m)})(y).
\end{eqnarray}
Thus it suffices to show that the character values $ \eta_n
(\lambda (\pi) ) (\alpha, s)$ and $\lambda (\C^{kn}) (\alpha, s)$
are equal for $(\alpha, s) \in \Gn,$ where $\alpha =(g, 1, \ldots,
1) \in \G^n$ and $s$ is an $n$-cycle, say $s = (12 \ldots n)$. It
is known (cf. \cite{FJW}) that the character value of $\eta_n
(\wt)(\alpha, s)$ is $\wt ( c)$, where $\wt$ is a class function
of $\G$ and $c$ is the conjugacy class of $g$. In particular
\begin{eqnarray*}
\eta_n (\lambda (\pi) )(\alpha, s) = \lambda (\pi) (c).
\end{eqnarray*}

Denote by $x^i_a, a =1, \ldots, n, i = 1, \ldots, k$ the
coordinates in $\C^{kn} = (\C^k)^n.$ The $a$-th factor $\G$ in
$\G^n \subset \Gn$ acts on $\C^k$ with coordinates $x^i_a, i =1,
\ldots, k$.

Consider the exterior monomial basis for $\Lambda^i \C^{kn}$.
Given such a monomial $X$, we first observe that the coefficient
of $X$ in $(\alpha, s) .X$ is $0$ unless there are equal numbers
of lower subscripts $1, 2, \ldots, n$ for $x^i_a$ appearing in
$X$. It follows that
\begin{eqnarray}  \label{eq_vanish}
\Lambda^i \C^{kn} (\alpha, s) = \mbox{trace }(\alpha, s)
\mid_{\Lambda^i \C^{kn}} = 0, \quad \mbox{if } i \mbox{ is not
divisible by } n.
\end{eqnarray}

For $i =m n, 1\leq m \leq k$, we further observe that the
coefficient of $X$ in $(\alpha, s) .X$ is $0$ unless the monomial
$X$ is of the form
\begin{eqnarray*}
 X(i_1 \ldots i_m) & =& x_1^{i_1} \wedge x_2^{i_1} \wedge \ldots
   \wedge x_n^{i_1} \wedge \\
   & & x_1^{i_2} \wedge x_2^{i_2} \wedge \ldots
   \wedge x_n^{i_2} \wedge \ldots \wedge x_1^{i_m} \wedge \ldots
   \wedge x_n^{i_m},
\end{eqnarray*}
where $\{ i_1, \ldots, i_m \}$ is an unordered $m$-tuple of
distinct numbers among $1, 2, \ldots, k$. Write $ \pi( g) x^{i}_1
= \sum_j b_{ij} x^j_1 $ and denote by $B$ the $k \times k$ matrix
$(b_{ij})$. Then
\begin{eqnarray*}
 (\alpha, s). X(i_1 \ldots i_m)
 & =& x_2^{i_1} \wedge x_3^{i_1} \wedge \ldots
   \wedge x_n^{i_1} \wedge \sum_j b_{i_1 j} x_1^j \wedge \\
   & &  x_2^{i_2} \wedge x_3^{i_2} \wedge \ldots
   \wedge x_n^{i_2} \wedge \sum_j b_{i_2 j} x_1^j
 \wedge \ldots \wedge \\
 & & x_2^{i_m} \wedge x_3^{i_m} \wedge \ldots
   \wedge x_n^{i_m} \wedge \sum_j b_{i_m j} x_1^j .
\end{eqnarray*}

It follows that the coefficient of $X(i_1 \ldots i_m)$ in
$(\alpha, s). X(i_1 \ldots i_m)$ is equal to
\begin{eqnarray}  \label{eq_det}
 && \sum_{\sigma} (-1)^{(n -1) m}
(-1)^{l (\sigma)} b_{i_1 \sigma (i_1)}b_{i_2 \sigma (i_2)} \ldots
b_{i_r \sigma (i_m)}    \nonumber \\
 & =&  (-1)^{(n -1) m} \det B(i_1 \ldots i_m),
\end{eqnarray}
where the summation runs over all permutations $\sigma$ of $i_1,
\ldots, i_m$, $l (\sigma) $ is the length of $\sigma$, and $\det
B(i_1 \ldots i_m)$ denotes the determinant of the $m \times m$
minor of $A$ consisting of the rows and columns $i_1, \ldots,
i_m$.

By (\ref{eq_vanish}) and (\ref{eq_det}) we calculate that
\begin{eqnarray*}   \label{eq_two}
\lambda (\C^{kn}) (\alpha, s)
 & =& \sum_{m =0}^k (-1)^{mn}\Lambda^{mn} (\C^{kn})(\alpha, s) \nonumber \\
 & =& \sum_{m =0}^k (-1)^{mn} (-1)^{m(n -1)} \sum_{ \{i_1, \ldots , i_m\} }
  \det B(i_1 \ldots i_m) \nonumber \\
 & =& \sum_{m =0}^k (-1)^m e_m (t_1, \ldots, t_k) \nonumber   \\
 & =& \sum_{m =0}^k (-1)^m (\Lambda^m \pi) (g) \nonumber   \\
 & =& \lambda (\pi) (c ),
 \end{eqnarray*}
where the summation runs over all unordered $r$-tuples $ \{i_1,
\ldots , i_m\}$ of distinct numbers among $1, 2, \ldots, k$, and
$e_m$ denotes the $m$-th elementary symmetric polynomial of the
eigenvalues $t_1, \ldots, t_k$ of the matrix $\pi (g) \in SL_k
(\C)$. The two identities involving $e_r$ used above are well
known.

By comparing the character values $\eta_n (\lambda (\pi) )(\alpha,
s)$ and $\lambda (\C^{kn})(\alpha, s) $ calculated above, we see
that
\[
\eta_n (\lambda (\pi) )(\alpha, s) =\lambda (\C^{kn})(\alpha, s) =
\lambda (\pi) (c ).
\]

% Since the sign character of $\Gn$ takes value
% $(-1)^{n -1}$ at $c_n$, we have
% $$
%  \varepsilon_n (\g) (c_n) = ( -1)^{n -1} \g (c).
% $$

Therefore if $x \in \Gn$ is of type $\rho \in {\cal P}_n (\G_*)$
then by using (\ref{eq_mult}) we obtain that
 \begin{eqnarray*}
  \eta_n (\lambda (\pi)  ) ( x)
    & =& \prod_{c\in \G_*} \lambda (\pi) (c)^{l (\rho(c))} \label{eq_term}\\
    & =& \lambda (\C^{kn}) (x).
%  \varepsilon_n (\g ) ( x)
%    = (-1)^n  \prod_{c\in \G_*} ( - \g (c))^{l (\rho(c))},
%   \label{eq_signterm}
 \end{eqnarray*}
 This completes the proof.
\end{demo}

\begin{remark}  \rm
The identification of the two virtual characters can be also seen
alternatively as follows:
\begin{eqnarray}
\la (\C^{kn}) & =& \sum_{i =0}^{kn} (-1)^i \Lambda ( (\C^k)^n )
 \nonumber   \\
 & =&  \sum_{i_1, \ldots , i_n } (-1)^{i_1 + \ldots + i_n}
 \Lambda^{i_1} (\C^k) \otimes \ldots \otimes \Lambda^{i_n} (\C^k)
  \nonumber  \\
 & =& \sum_{\{n_0, \ldots , n_k \}} (-1)^{\sum_i i n_i}
  \mbox{Ind}^{\Gn}_{\G_{n_0} \times \ldots \times \G_{n_k}}
  \Lambda^0 (\C^k)^{\boxtimes n_0} \otimes \ldots \otimes
  \Lambda^k (\C^k)^{\boxtimes n_k} \label{eq_def}  \\
 & =& (\sum_i (-1)^i \Lambda^i (\C^k) )^{\boxtimes n}   \label{eq_work}\\
 & =& (\la (\C^k) )^{\boxtimes n}  \nonumber
\end{eqnarray}
where $\{n_0, \ldots , n_k \}$ ranges over the $(k+1)$-tuple of
non-negative integers such that $ \sum_i n_i = n.$
Eq.~(\ref{eq_def}) above basically follows from the definition of
an induction functor. Eq.~(\ref{eq_work}) is a generalization of
(\ref{eq_virt}) which can be established with some effort.
\end{remark}

When $\G$ is trivial then $\la (\pi) = 0 \in R(\G)$ and so $\eta_n
(\la (\pi)) =0$. We have an immediate corollary.
\begin{corollary}
When $\G$ is trivial and $\Gn$ becomes the symmetric group $S_n$,
the virtual $S_n$-character $\lambda (\C^{kn})$ is zero.
\end{corollary}
\subsection{Derived categories and Grothendieck Groups}
In this subsection we let $\G $ be a finite subgroup of $SL_2
(\C)$ unless otherwise specified.

We denote by $D_{\Gn} (\C^{2n})$ the bounded derived category of
$\Gn$-equivariant coherent sheaves on $\C^{2n}$, and denote by
$D(\alen)$ the bounded derived category of coherent sheaves on
$\alen$. Define two functors $\Phi : D(\alen)  \rightarrow D_{\Gn}
(\C^{2n})$ and $\Psi : D_{\Gn} (\C^{2n}) \rightarrow D(\alen) $ by
\begin{eqnarray*}
 \Phi (-)  & =& Rp_* ({\cal O}_{\univ} \bigotimes q^*(-))\\
 \Psi (-)  & =& (Rq_* RHom({\cal O}_{\univ}, p^* (-)) )^{\Gn},
\end{eqnarray*}
where $\univ$ is the universal scheme defined in (\ref{eq_univ}),
and $p, q$ denote the projections $\alen \times \C^{2n}$ to
$\C^{2n}$ and $\alen$ respectively.

Take a basis of $R(\Gn)$ given by the irreducible characters
$s_{\rho}, \rho \in {\cal P}_n (\G^*)$ of $\Gn$ (cf. \cite{M, Z}).
We denote by ${\cal O}_{\C^{2n}}$ the structure sheaf over
$\C^{2n}$. Recall that the vector bundle (i.e. locally free sheaf)
${\cal V}_{\rho}$ is defined in (\ref{eq_wreathreg}). The
following theorem can be derived by using a general theorem due to
Bridgeland, King and Reid (Theorem~1.2, \cite{BKR}) since $\tau_n
: \ale^{[n]} \rightarrow \C^{2n} /\Gn$ is a crepant resolution and
$\Gn$ preserves the symplectic structure of $\C^{2n}$.

\begin{theorem}
Under the assumption of the validity of
Conjecture~\ref{conj_main}, $\Phi$ is an equivalence of categories
and $\Psi$ is its adjoint functor. In particular, $\Psi$ sends
${\cal O}_{\C^{2n}} \otimes s_{\rho}^{\vee}$ to ${\cal V}_{\rho}$.
\end{theorem}

\begin{remark} \rm
When $n=1 $, Conjecture~\ref{conj_main} is known to be true (cf.
Remark~\ref{rem_one}) and the above theorem was established by
Kapranov and Vasserot \cite{KV}.
\end{remark}

\begin{remark} \rm \label{rem_alt}
By assuming the validity of the $n!$ conjecture (cf. \cite{GH}),
we have by Remark~\ref{rem_one} that $\ale^n//S_n \cong \alen$.
Thus we may replace $\alen$ by $\ale^n // S_n$ in the crepant
Hilbert-Chow resolution $\alen \stackrel{\pi_n}{\rightarrow}
\ale^n / S_n$. Clearly $S_n$ preserves the holomorphic symplectic
structure of the Cartisian product $\ale^n$. Then one can apply
again Theorem~1.2 in \cite{BKR} to show that there is an
equivalence of derived categories between $D(\alen)$ and $D_{S_n}
(\ale^n)$ of $S_n$-equivariant coherent sheaves of $\alen$.
\end{remark}

Below we further specialize and apply the general results of
Bridgeland, King and Reid \cite{BKR} to our setup. We denote by
$D^0_{\Gn} (\C^{2n})$ the full subcategory of $D_{\Gn} (\C^{2n})$
consisting of objects whose cohomology sheaves are concentrated on
the origin of $\C^{2n}$. We denote by $D^0_{S_n} (\ale^n)$ the
full subcategory of $D_{S_n} (\ale^n)$ consisting of objects whose
cohomology sheaves are concentrated on the $n$-th Cartisian
product of the exceptional divisor $D \subset \ale$. Denote by
$D^0(\alen)$ the full subcategory of $D(\alen)$ consisting of
objects whose cohomology sheaves are concentrated on the fiber
$\ale^{[n], 0}$ of $\alen$. Recall that $\ale^{[n], 0}$ is
described in Subsection~\ref{subsect_fiber}.

The equivalence between the derived categories $D_{\Gn} (\C^{2n})$
and $D(\ale^{[n]})$ induces an equivalence between $D^0_{\Gn}
(\C^{2n})$ and $D^0(\ale^{[n]})$. Thus we have the following
commutative diagram (under the assumption of the validity of
Conjecture~\ref{conj_main}):
\begin{eqnarray}   \label{eq_cat}
 \begin{array}{ccc}
   D^0_{\Gn}(\C^{2n})   & \stackrel{\simeq}{\longrightarrow}
      & D^0(\ale^{[n]})   \\
   \downarrow  &                 & \downarrow   \\
   D_{\Gn}(\C^{2n})     & \stackrel{\simeq}{\longrightarrow}
      & D(\ale^{[n]}) .
 \end{array}
\end{eqnarray}
Given objects $E, F$ in $D_{\Gn}(\C^{2n})$ and
$D^0_{\Gn}(\C^{2n})$ respectively, we define the Euler
characteristic
\begin{eqnarray*}
\chi^{\Gn} (E, F) = \sum_i ( -1)^i \dim Hom_{D_{\Gn}(\C^{2n})} (E,
F[i]).
\end{eqnarray*}
This gives a natural bilinear pairing between $D_{\Gn} (\C^{2n})$
and $D^0_{\Gn} (\C^{2n})$. Similarly we can define the Euler
characteristic $\chi (A, B)$ for objects $A, B$ in $D(\ale^{[n]})$
and $D^0(\ale^{[n]})$ respectively, which gives rise to a bilinear
pairing between $D(\ale^{[n]})$ and $D^0(\ale^{[n]})$. We further
have $\chi^{\Gn}(E, F) =\chi (\Psi (E), \Psi (F)),$ cf.
\cite{BKR}.

We denote by $K_{\Gn} (\C^{2n}), K^0_{\Gn} (\C^{2n}), K
(\ale^{[n]})$, $ K^0 (\ale^{[n]})$, $K^0_{S_n} (\ale^n)$ and
$K_{S_n} (\ale^n)$ the Grothendieck groups of the corresponding
derived categories. It is well known that $K_{\Gn}
(\C^{2n})$ and $ K^0_{\Gn} (\C^{2n})$ are both isomorphic to the
representation ring $R_{\Z}( \Gn)$. The bilinear pairings
mentioned above together with the embeddings of categories induces
a bilinear form on $K^0_{\Gn} (\C^{2n})$ and respectively on
$K^0(\ale^{[n]})$.

Let ${\cal O}_0$ be the skyscraper sheaf at the origin $0$ on
$\C^{2n}$. The $\Gn$-bundles ${\cal O} \otimes s_{\rho}, \rho \in
{\cal P}_n (\G^*)$ form a basis for $K_{\Gn} (\C^{2n})$ while the
modules $s_{\rho} \otimes {\cal O}_0, \rho \in {\cal P}_n (\G^*)$
form the dual basis for $K^0_{\Gn} (\C^{2n})$.
\begin{theorem} \label{th_reform}
The map by sending $s_{\rho}$ to ${\cal O}_0 \otimes s_{\rho}$,
$\rho \in {\cal P}_n (\G^*)$, is an isometry between $R_{\Z}(\Gn)$
endowed with the weighted bilinear form and $K^0_{\Gn} (\C^{2n})$
endowed with the bilinear form defined above. In particular the
bilinear form on $K^0_{\Gn} (\C^{2n})$ is semipositive definite
symmetric.
\end{theorem}

\begin{demo}{Proof}
Following a similar argument as in Gonzalez-Sprinberg and Verdier
\cite{GSV} (which is for $n =1$), we obtain the following
commutative diagram by using the Koszul resolution of ${\cal O}_0$
on $\C^{2n}$:
\begin{eqnarray*}
 \begin{array}{ccc}
   R_{\Z}(\Gn)    & \stackrel{\simeq}{\longrightarrow}
     &  K^0_{\Gn} (\C^{2n} )     \\
     \jmath \downarrow  &                 & \downarrow  \\
    R_{\Z}(\Gn)    & \stackrel{\simeq}{\longrightarrow}
     &  K_{\Gn} (\C^{2n} ).
 \end{array}
\end{eqnarray*}
Here the horizontal maps are isomorphisms given by sending
$s_{\rho}$ to ${\cal O}_0 \otimes s_{\rho}$ and respectively to
${\cal O}_{\C^{2n}} \otimes s_{\rho}$, $\rho \in {\cal P}_n
(\G^*)$. The left vertical map $ \jmath $ is given by
multiplication by the virtual character $\lambda (\C^{2n})$ of
$\Gn$, and the right vertical one is induced from the natural
embedding of the corresponding categories. Now the theorem follows
from the definition of the weighted bilinear form on $R(\Gn)$,
Theorem~\ref{th_key}, and the fact that the basis ${\cal O}_0
\otimes s_{\rho}, \rho \in {\cal P}_n (\G^*)$ for $K^0_{\Gn}
(\C^{2n})$ is dual to the basis ${\cal O}_{\C^{2n}} \otimes
s_{\rho}, \rho \in {\cal P}_n (\G^*)$ for $K_{\Gn} (\C^{2n})$.
\end{demo}

\begin{remark} \rm
The main results in \cite{FJW} (see Theorem 7.2 and Theorem 7.3 in
{\em loc. cit.}) can be now formulated by using the space
\[
{\cal F}_{\G} = \bigoplus_{n \geq 0}  K^0_{\Gn} (\C^{2n})
\bigotimes \C [K^0_{\G}(\C^2 )]
\]
with its natural bilinear form induced from the Koszul-Thom class.
Here $ \C [ - ]$ denotes the group algebra. Roughly speaking,
${\cal F}_{\G}$ affords a vertex representation of the toroidal
Lie algebra and a distinguished subspace of ${\cal F}_{\G}$
affords the basic representation of the affine Lie algebra
$\affineg$ whose associated affine Dynkin diagram corresponds to
$\G$ in the sense of McKay. This may be viewed as a form of McKay
correspondence relating finite subgroups of $SL_2 (\C)$ to affine
and toroidal Lie algebras.
\end{remark}

Now we have the following commutative diagram (assuming the
validity of Conjecture~\ref{conj_main}):
\begin{eqnarray*}
 \begin{array}{ccccc}
   R_{\Z}(\Gn)    & \stackrel{\simeq}{\longrightarrow}
     &  K^0_{\Gn} (\C^{2n} )   & \stackrel{\simeq}{\longrightarrow}
     &   K^0 (\ale^{[n]} ) \\
     \jmath \downarrow  &   & \downarrow &     & \downarrow  \\
    R_{\Z}(\Gn)    & \stackrel{\simeq}{\longrightarrow}
     &  K_{\Gn} (\C^{2n} ) & \stackrel{\simeq}{\longrightarrow}
     & K (\ale^{[n]} ).
 \end{array}
\end{eqnarray*}
By Remark~\ref{rem_alt} and a similar argument as above, we have
another commutative diagram (assuming the validity of
Conjecture~\ref{conj_main}):
\begin{eqnarray*}
 \begin{array}{ccc}
     K^0_{S_n} (\ale^n )   & \stackrel{\simeq}{\longrightarrow}
     &   K^0 (\ale^{[n]} ) \\
      \downarrow &     & \downarrow  \\
      K_{S_n} (\ale^n) & \stackrel{\simeq}{\longrightarrow}
     & K (\ale^{[n]} ).
 \end{array}
\end{eqnarray*}

Combining the two diagrams above, we have obtained the following
theorem.

\begin{theorem}
Under the assumption of the validity of
Conjecture~\ref{conj_main}, the isomorphisms among $(R_{\Z}(\Gn),
\langle -, - \rangle_{\lambda(\C^2)} )$, the K-groups $K^0_{S_n}
(\ale^{n} )$ and $K^0(\ale^{[n]})$  are  isometries.
\end{theorem}

\begin{remark}  \rm
All the algebraic structures on $\bigoplus_n R (\Gn)$ (cf.
\cite{W,FJW}) and thus on $\bigoplus_n K^0_{\Gn} (\C^{2n})$ can
now be carried over to $\bigoplus_n K^0(\ale^{[n]})$. However it
remains to match these with the Grojnowski-Nakajima construction
on $\bigoplus_{n} H(\alen)$ in terms of correspondence varieties.
\end{remark}

\begin{remark}  \rm
Given a finite subgroup $G$ of $SL_K (\C)$, one asks whether there
is a crepant resolution $Y$ of the affine orbifold $\C^K / G$ so
that there exists a canonical isomorphism between $K_{G}(\C^K) =
K(Y)$; one further asks whether or not the answer can be provided
by the Hilbert quotient $\C^K // G$, cf. Reid \cite{R}. The answer
is affirmative for $K =2$, known as the McKay correspondence
\cite{GSV} (also compare \cite{IN, KV}). For $K=3$, there has been
much work by various people, cf. \cite{Ro, R, Nr, IN} and
references therein, and it is settled by Bridgeland-King-Reid
\cite{BKR}. However not much is known in general (see however
\cite{BKR}) and there has been counterexamples. Our work provides
strong evidence for an affirmative answer in the case of $\C^{2n}$
acted upon by the wreath product $\Gn$ associated to a finite
subgroup $\G \subset SL_2 (\C)$.
\end{remark}
\section{A direct isomorphism of algebraic structures on equivariant K-theory}
\label{sect_ktheory}
In this section we assume that the reader is familiar with
\cite{W}. For shortness of notations we will use $K^{\topo}(-)$
and $K^{\topo}_G(-)$ to denote the {\em complexified}
($G$-equivariant) topological K-group. We further assume that $X$
is a quasi-projective surface acted upon by a finite subgroup $\G$
and $Y$ is a resolution of singularities of $X/ \G$ such that
there exists a canonical isomorphism $\theta$ between
$K^{\topo}_{\G}(X)$ and $K^{\topo}( Y)$.

Let ${\cal C} = \{ V_1, \ldots, V_l \}$ be a basis for $\Kgx$.
Without loss of generality we may assume they are genuine
$\G$-vector bundles on $X$. We denote by $W_i = \theta ( V_i)$.
The set $ \{ W_1, \ldots, W_l \}$ is a basis for $\Ky$. We remark
that representatives of $W_i$'s can again be chosen as certain
canonical vector bundles in favorable cases, including the
important case when $X= \C^2$ and $\G \subset SL_2 (\C )$.

Let $S_{\la} $ be the irreducible representation of the symmetric
group $S_n $ associated to the partition $\lambda$ of $n$. Define
$$
  S_{\la}  (V_i) = S_{\la}  \bigotimes V_i^{ \boxtimes n}.
$$ Endowed with the diagonal action of $S_n$ on the two tensor
factors and the action of $\G^n$ on the second factor, $S_{\la}
(V_i)$ is an $\Gn$-equivariant vector bundle.

Given a partition-valued function $\la = (\la_i )_{ 1 \leq i \leq
l} \in {\cal P}_n (\cal C)$, we define the $\Gn$-equivariant
vector bundle $$ S^X_{\la} =
 \Ind^{\Gn}_{ \G_{|\la_1|} \times \ldots \times \G_{|\la_rl|} }
  S_{\la_1}(V_1 ) \times \ldots \times S_{\la_l}(V_l).
$$
In a parallel way, we can define the $S_n$-equivariant bundle $
S^Y_{\la}$ associated to $ \la \in {\cal P}_n (\cal C) $ as $$
S^X_{\la}
=
 \Ind^{S_n}_{ S_{|\la_1|} \times \ldots \times S_{|\la_l|} }
  S_{\la_1}(\theta(V_1) ) \times \ldots \times S_{\la_l}(\theta (V_l)).
$$

Recall that we constructed in \cite{W} various algebraic
structures such as Hopf algebra, $\lambda$-ring, Heisenberg
algebra on $\bigoplus_{n \geq 0} \Kgxn$ for a $\G$-space $X$,
generalizing the results of Segal \cite{S2} for $\G$ trivial. The
following proposition can be proved using \cite{W} in a way as
Macdonald \cite{M} did when $X$ is a point.

\begin{proposition}
The $\Gn$-bundles $ S^X_{\la}, \la \in {\cal P}_n (\cal C)$ form a
basis of $\Kgxn$. The $S_n$-bundles $ S_{\la}^Y, \la \in  {\cal
P}_n (\cal C)$ form a basis of $\Kyn$.
\end{proposition}
These bases will be referred to as Schur bases, generalizing the
usual one for $R(\Gn) = K^{\topo}_{\Gn} (pt)$.

\begin{theorem}
The map $\Theta$ from $\bigoplus_{n \geq 0} \Kgxn$ to
$\bigoplus_{n \geq 0} \Kyn$ by sending $S_{\la} ^X$ to $S_{\la}
^Y$ is an isomorphism of Hopf algebras, $\la$-rings, and
representations over the Heisenberg algebra.
\end{theorem}

\begin{demo}{Proof}
We use \cite{W} as a basic reference. We follow the notations
there with an additional use of $X, Y$ as subscripts to specify
the space we are referring to.

As graded vector spaces $\bigoplus_{n \geq 0} \Kgxn$ and
$\bigoplus_{n \geq 0} \Kyn$ have the same graded dimension due to
the isomorphism between $K^{\topo}_{\G} (X)$ and $K^{\topo} (Y)$,
cf. Theorem 3 in \cite{W}. Thus the map $\Theta$ given by matching
the Schur basis is an additive isomorphism.

Recall that the Adam's operations $\varphi^m_X$ on the space
$\bigoplus_{n \geq 0} \Kgxn$ and $\varphi^m_Y$ on $\bigoplus_{n \geq
0} \Kyn$ satisfy the identities (where $q$ is a formal parameter),
cf. Proposition 4 in \cite{W}:
\begin{eqnarray*}
\bigoplus_{n \geq 0} q^n V_i^{\boxtimes n}
  & =& \exp \left( \sum_{m > 0} \frac1m  \varphi^m_X (V_i) q^m
            \right),   \\
\bigoplus_{n \geq 0} q^n W_i^{\boxtimes n}
  & =& \exp \left( \sum_{m > 0} \frac1m \varphi^m_Y (W_i) q^m
            \right).
\end{eqnarray*}

It follows that $\varphi^m_X (V_i)$ and $\varphi^m_Y (W_i)$ are
uniquely determined by $V_i^{\boxtimes n}$ $(n \geq 0)$ and
respectively $W_i^{\boxtimes n}$ in the same way (by taking
logarithms of the above identities). Since the isomorphism
$\Theta$ sends $V_i^{\boxtimes n}$ to $W_i^{\boxtimes n}$, the
Adams operations $\phi^r_X (V_i)$ and $\phi^r_X (V_i)$ matches
under $\Theta$ and so does the $\lambda$-ring structures on
$\bigoplus_{n \geq 0} \Kgxn$ and $\bigoplus_{n \geq 0} \Kyn$.

Recall that Heisenberg algebras ${\cal H}_X$ and ${\cal H}_Y$ constructed
in terms of K-theory maps act
irreducibly on $\bigoplus_{n \geq 0} \Kgxn$ and $\bigoplus_{n \geq
0} \Kyn$ respectively, cf. Theorem 4 in \cite{W}. The Heisenberg
algebra generators are essentially defined in terms of Adams operations,
induction functors and restriction functors such as
$\mbox{Ind}_{\G_m \times \G_{n -m}}^{\Gn}$, $\mbox{Res}_{\G_m
\times \G_{n -m}}^{\G_n}$, etc. Since the Adams operations,
induction and restriction functors are compatible with
the Schur bases and thus with the map
$\Theta$, the Heisenberg algebras acting on $\bigoplus_{n \geq 0}
\Kgxn$ and $\bigoplus_{n \geq 0} \Kyn$ also matches under
$\Theta$.
\end{demo}

Department of Mathematics, North Carolina State Univer\-sity,
Raleigh, NC 27695-8205. wqwang@math.ncsu.edu


\frenchspacing
\begin{thebibliography}{FLM}

\bibitem[BG]{BG} R. Bezrukavnikov and V. Ginzburg,
{\em Hilbert schemes and reductive groups}, in preparation.

\bibitem[BKR]{BKR} T. Bridgeland, A. King and M. Reid,
{\em Mukai implies McKay}, preprint, math.AG/9908027.

\bibitem[CM]{CM} M. de Cataldo and L. Migliorini,
{\em The Douady space of a complex surface},
math.AG/9811159, to appear in Adv. in Math.

\bibitem[Fo]{Fo} J. Fogarty,
{\em Algebraic families on an algebraic surface}, Amer. J. Math.
{\bf 90} (1968) 511--521.

\bibitem[FJW]{FJW} I.~B.~Frenkel, N.~Jing
and W.~Wang, {\em Vertex representations via finite groups and the
McKay correspondence}, math.QA/9907166, to appear in IMRN.

\bibitem[FK]{FK} I.~B. Frenkel and V.~G. Kac, {\em
Basic representations of affine Lie algebras and dual resonance
models}, Invent. Math. {\bf 62} (1980) 23--66.

\bibitem[GH]{GH} A. Garsia and M. Haiman,
{\em A graded representation model for Macdonald's polynomials},
Proc. Nat. Acad. USA {\bf 903} (1993) 3607--3610.

\bibitem[GSV]{GSV} G.~Gonzalez-Sprinberg and J.-L.~Verdier,
{\em Construction g\'eom\'etrique de la correspondance de McKay},
Ann. Sci. \'Ecole Norm. Sup. {\bf 16} (1983) 409--449.

\bibitem[G]{G} L. G\"ottsche,
{\em The Betti numbers of the Hilbert scheme of points on a smooth
projective surface}, Math. Ann. {\bf 286} (1990) 193--207.

\bibitem[Gr]{Gr} I.~Grojnowski,
{\em Instantons and affine algebras I: the Hilbert scheme and
vertex operators}, Math. Res. Lett. {\bf 3} (1996) 275--291.

\bibitem[H]{H} M. Haiman,
{\em Macdonald polynomials and Hilbert schemes}, UCSD preprint.

\bibitem[IN]{IN} Y. Ito and H. Nakajima,
{\em McKay correspondence and Hilbert schemes in dimension three},
math.AG/9803120, to appear in Topology.

\bibitem[INr]{INr} Y. Ito and I. Nakamura,
{\em McKay correspondence and Hilbert schemes}, Proc. Japan Acad.
Ser. {\bf  A 72} (1996) 135--138.

\bibitem[Ka]{Ka} M. Kapranov, {\em Chow quotients of
Grassmannians}, in Gelfand Seminar {\bf 1} (eds. S.~Gelfand,
S.~Gindikin), Adv. in Soviet Math. {\bf 16} (1993) 29--110, Amer.
Math. Soc. Providence RI.

\bibitem[KV]{KV} M. Kapranov and E.~Vasserot,
{\em Kleinian singularities, derived categories and Hall
algebras}, preprint, math.AG/9812016.

\bibitem[Kr]{Kr} P.~Kronheimer,
{\em The construction of ALE spaces as hyper-K\"ahler quotients},
J. Diff. Geom. {\bf 28} (1989) 665--683.

\bibitem[KN]{KN} P.~Kronheimer and H.~Nakajima,
{\em Yang-Mills instantons on ALE gravitational instantons}, Math.
Ann. {\bf 288} (1990) 263--307.

%\bibitem[L]{L} M. Lehn,
%{\em Chern classes of tautological sheaves on Hilbert schemes of
%points on surfaces}, Invent. Math. {\bf 136} (1999) 157--207.

\bibitem[M]{M} I.~G. Macdonald,
{\em Polynomial functors and wreath products}, J. Pure Appl. Alg.
{\bf 18} (1980) 173--204.

\bibitem[Mc]{Mc} J.~McKay, {\em Graphs, singularities and finite
groups}, Proc. Sympos. Pure Math. {\bf 37}, Amer. Math. Soc,
Providence, RI (1980) 183--186.

\bibitem[N]{N} H. Nakajima,
{\em Instantons on ALE spaces, quiver varieties, and Kac-Moody
algebras}, Duke Math. J. {\bf 76} (1994) 365--416.

\bibitem[N1]{N1} H. Nakajima,
{\em Heisenberg algebra and Hilbert schemes of points on
projective surfaces}, Ann. Math. {\bf 145} (1997) 379--388.

\bibitem[N2]{N2} H. Nakajima,
{\em Lectures on Hilbert schemes of points on surfaces}, to be
published by AMS.

\bibitem[N3]{N3} H. Nakajima, {\em Quiver varieties and Kac-Moody
algebras}, Duke Math. J. {\bf 91} (1998) 515--560.

\bibitem[N4]{N4} H. Nakajima, {\em Private communications}.

\bibitem[Nr]{Nr} I. Nakamura,
{\em Hilbert schemes of abelian group orbits}, to appear in J.
Alg. Geom.

\bibitem[R]{R} M. Reid,
{\em McKay correspondence}, in Proc. alg. geom sympos. (Kinosaki,
Nov. 1996), T. Katsura (ed.), 14-41; alg-geom/9702016.

\bibitem[Ro]{Ro} S.-S. Roan,
{\em Minimal resolutions of Gorenstein orbifolds in dimension
three}. Topology 35 (1996) 489--508.

\bibitem[S1]{S1} G.~Segal, {\em Unitary representations of
some infinite dimensional groups},
Commun. Math. Phys. {\bf 80} (1981) 301--342.

\bibitem[S2]{S2} G.~Segal, {\em Equivariant K-theory and
symmetric products}, 1996 preprint (unpublished).

\bibitem[VV]{VV} M. Varagnolo and E. Vasserot,
{\em On the K-theory of the cyclic quiver variety}, preprint,
math.AG/9902091.

\bibitem[VW]{VW} C. Vafa and E. Witten,
{\em A strong coupling test of $S$-duality}, Nucl. Phys. {\bf B
431} (1994) 3--77.

\bibitem[W]{W} W.~Wang, {\em Equivariant K-theory and wreath products},
MPI preprint {\bf \#~86}, August 1998; {\em Equivariant K-theory,
wreath products and Heisenberg algebra}, math.QA/9907151, to
appear in Duke Math. J.

\bibitem[Wey]{Wey} H. Weyl, {\em The classical groups, their
invariants and representations}, Princeton University Press, 1946.

\bibitem[Z]{Z} A.~Zelevinsky, {\em Representations of finite classical
groups. A Hopf algebra approach}, Lect. Notes in Math. {\bf 869},
Springer-Verlag, Berlin-New York, 1981.

\end{thebibliography}
\end{document}